\documentclass[12pt]{article}
\usepackage{amssymb, amsmath, amsfonts, tikz, geometry, ytableau, delimset}
\usepackage[indent,margin=1cm]{caption}
\usepackage[colorlinks,linkcolor=blue,anchorcolor=blue,citecolor=blue]{hyperref}
\numberwithin{equation}{section}
\numberwithin{figure}{section}

\hypersetup{colorlinks = true}
\newtheorem{thm}{Theorem}[section]
\newtheorem{conj}[thm]{Conjecture}
\newtheorem{conjc}[thm]{Conjollary}

\newtheorem{lem}[thm]{Lemma}
\newtheorem{prop}[thm]{Proposition}

\def\pf{\noindent{\it Proof.} }
\def\qed{\nopagebreak\hfill{\rule{4pt}{7pt}}
\medbreak}

\allowdisplaybreaks

\usepackage{tikz, tkz-berge}

\tikzset{OuterBoundary/.style={very thick}}
\tikzset{RimHook/.style={cyan, line width=1.2pt, rounded corners=2pt}}
\tikzset{edge/.style={cyan, very thick}}
\tikzset{ball/.style={shape=circle, inner sep=9pt, ball color=magenta!25!white}}
\tikzset{RHiball/.style={ball color=magenta!50!white}}
\tikzset{RHeball/.style={ball color=magenta!25!white}}
\tikzset{apball/.style={ball color=black}}

\SetVertexNoLabel
\tikzstyle{EdgeStyle}=[very thick, color=cyan]
\tikzstyle{VertexStyle}=[shape=circle, inner sep=6pt, ball color=magenta!25!white]

\usepackage{genyoungtabtikz}
\YFrench
\Yboxdim{1cm}

\begin{document}

\begin{center}
{\large \bf The twinning operation on graphs does not always preserve $e$-positivity}
\end{center}

\begin{center}
Ethan Y.H. Li$^{1}$, Grace M.X. Li$^{2}$, David G.L. Wang$^3$ and Arthur L.B. Yang$^{4}$\\[6pt]

$^{1,2,4}$Center for Combinatorics, LPMC\\
Nankai University, Tianjin 300071, P. R. China\\[8pt]

$^{3}$School of Mathematics and Statistics, Beijing Key Laboratory on MCAACI\\
Beijing Institute of Technology, Beijing 102488, P. R. China\\[8pt]

Email: $^{1}${\tt yinhao\_li@mail.nankai.edu.cn}, $^{2}${\tt limengxing@mail.nankai.edu.cn}, $^3${\tt glw@bit.edu.cn}, $^{4}${\tt yang@nankai.edu.cn}
\end{center}

\noindent\textbf{Abstract.}
Motivated by Stanley's $\mathbf{(3+1)}$-free conjecture on chromatic symmetric functions, Foley, Ho\`{a}ng and Merkel introduced the concept of strong $e$-positivity and conjectured that a graph is strongly $e$-positive if and only if it is (claw, net)-free. In order to study strongly $e$-positive graphs, they further introduced the twinning operation on a graph $G$ with respect to a vertex $v$, which adds a vertex $v'$ to $G$ such that $v$ and $v'$ are adjacent and any other vertex is adjacent to both of them or neither of them. Foley, Ho\`{a}ng and Merkel conjectured that if $G$ is $e$-positive, then so is the resulting twin graph $G_v$ for any vertex $v$. Based on the theory of chromatic symmetric functions in non-commuting variables developed by Gebhard and Sagan, we establish the $e$-positivity of a class of graphs called tadpole graphs. By considering the twinning operation on a subclass of these graphs with respect to certain vertices we disprove the latter conjecture of Foley, Ho\`{a}ng and Merkel. We further show that if $G$ is $e$-positive, the twin graph $G_v$ and more generally the clan graphs $G^{(k)}_v$ ($k \ge 1$) may not even be $s$-positive, where $G^{(k)}_v$ is obtained from $G$ by applying $k$ twinning operations to $v$.



\section{Introduction}

The main objective of this paper is to provide counterexamples to a conjecture due to Foley, Ho\`{a}ng and Merkel \cite{FHM19}, which arose in their study of Stanley's celebrated $\mathbf{(3+1)}$-free conjecture on chromatic symmetric functions. Their conjecture aims to show the twinning operation on graphs preserves the positivity of the coefficients in the expansion of their chromatic symmetric functions in terms of elementary symmetric functions. In this paper we use tadpole graphs to construct a family of counterexamples to Foley, Ho\`{a}ng and Merkel's conjecture.

Let us first give an overview of some background. Let $\mathbb{P}$ be the set of positive integers. Let $G$ be a simple graph with vertex set $V(G) = \{v_1,\ldots,v_d\}$ and edge set $E(G)$. A \textit{proper coloring} of $G$ is a map $\kappa: V(G) \to \mathbb{P}$ such that $\kappa(u) \neq \kappa(v)$ for any two adjacent vertices $u$ and $v$. It is known that the number of proper colorings of $G$ with $n$ colors is a polynomial of $n$, which is called the \textit{chromatic polynomial} of $G$ and denoted by $\chi_G(n)$. Given a set ${\bf{x}} = \{x_1,x_2,\ldots\}$ of countably infinite commuting indeterminates,  Stanley \cite{Sta95} studied the generating function
\[
    X_G({\bf{x}})= \sum_{\kappa} x_{\kappa(v_1)}x_{\kappa(v_2)}\cdots x_{\kappa(v_d)},
\]
where the sum ranges over all proper colorings of $G$. It is clear that $X_G({\bf{x}})$ is a symmetric function of degree $d$,
and equals $\chi_G(n)$ when we set $x_i=1$ for $1\leq i\leq n$ and $x_i=0$ otherwise. For this reason,
$X_G({\bf{x}})$ is called the \textit{chromatic symmetric function} of $G$. It is natural to consider the expansion of
$X_G({\bf{x}})$ in terms of various basis of the algebra of symmetric functions, such as the elementary symmetric functions $e_{\lambda}({\bf{x}})$ and the Schur functions $s_{\lambda}({\bf{x}})$, which are indexed by integer partitions $\lambda$.
From now on we shall omit the indeterminate set ${\bf{x}}$ if there is no confusion. We say that $G$ is $e$-positive (resp.\ $s$-positive) if all the coefficients are non-negative in the expansion $X_G=\sum_{\lambda}c_{\lambda}e_{\lambda}$ (resp.\ $X_G=\sum_{\lambda}d_{\lambda}s_{\lambda}$).
Stanley's $\mathbf{(3+1)}$-free conjecture states that the incomparability graph of any $\mathbf{(3+1)}$-free partially ordered set (poset) is $e$-positive, which is equivalent to a conjecture due to Stanley and Stembridge on immanants of the Jacobi-Trudi matrices \cite{SS93}.
Recall that a poset is said to be \textit{$\mathbf(a+b)$-free} if it contains no induced subposet isomorphic to the disjoint union of an $a$-element chain and a $b$-element chain. The \textit{incomparability graph} of $P$ is the graph obtained by taking the elements of $P$ as its vertices and letting two vertices be adjacent if and only if they are not comparable in $P$. Stanley's $\mathbf{(3+1)}$-free conjecture has received considerable attention, see for instance \cite{DFW17,FKKAMT18,Gas98,GP13,Hai95,HHT19,HP18,Lee08,LY19,SW16,Ska01} and references therein.

\begin{figure}[ht]
\centering
\begin{tikzpicture}[scale = 1.5]
    \fill (0,0) circle (0.3ex);
    \draw (0,0) -- (0,-1);
    \fill (0,-1) circle (0.3ex);
    \draw (0,0) -- (-1,-1);
    \fill (-1,-1) circle (0.3ex);
    \draw (0,0) -- (1,-1);
    \fill (1,-1) circle (0.3ex);
    
    \node[below] (w) at (0,-1.2) {claw};
    
    \fill (3,0.3) circle (0.3ex);
    \draw (3,0.3) -- (3,-0.2) -- (2.5,-1) -- (2,-1) -- (3.5,-1) -- (4,-1);
    \draw (3,-0.2) -- (3.5,-1);
    \fill (3,-0.2) circle (0.3ex);
    \fill (2.5,-1) circle (0.3ex);
    \fill (3.5,-1) circle (0.3ex);
    \fill (2,-1) circle (0.3ex);
    \fill (4,-1) circle (0.3ex);
    
    \node[below] (t) at (3,-1.2) {net};
    
    \fill (6,0.3) circle (0.3ex);
    \fill (5,0.3) circle (0.3ex);
    \draw (6,0.3) -- (6,-0.2) -- (5.5,-1) -- (5,-0.2) -- (5,0.3);
    \draw (6,-0.2) -- (5,-0.2);
    \fill (6,-0.2) circle (0.3ex);
    \fill (5.5,-1) circle (0.3ex);
    \fill (5,-0.2) circle (0.3ex);
    
    \node[below] (l) at (5.5,-1.2) {bull};
\end{tikzpicture}
    \caption{The claw, net and bull graphs.}\label{clawnetbull}
\end{figure}
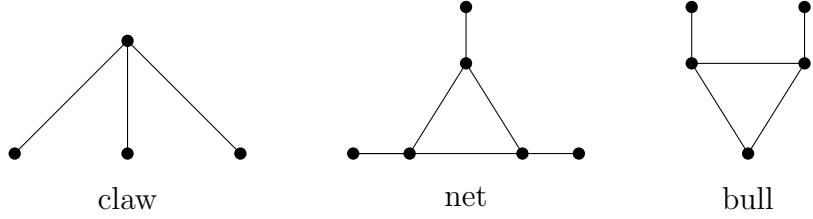

Stanley's $\mathbf{(3+1)}$-free conjecture also stirred up the interest of the combinatorial community to study the chromatic symmetric functions of claw-free graphs, since the claw, as depicted in Figure \ref{clawnetbull}, is the incomparability graph of the poset $\mathbf{(3+1)}$, and hence a graph $G$ is an incomparability graph of a $\mathbf{(3+1)}$-free poset if and only if it is a claw-free incomparability graph.
However, it is not true that every claw-free graph is $e$-positive, and one counterexample is the net graph given by Stanley \cite{Sta95}, as shown in Figure \ref{clawnetbull}. On account of this, Foley, Ho\`{a}ng and Merkel \cite{FHM19} initiated the study of (claw, net)-free graphs, namely, those graphs containing no induced subgraph isomorphic to the claw or the net. They conjectured that all (claw, net)-free graphs are $e$-positive. As implicitly stated in \cite{FHM19}, this conjecture is equivalent to saying that a graph is (claw, net)-free if and only if it is strongly $e$-positive, where a graph is called \textit{strongly $e$-positive} if all its induced subgraphs are $e$-positive. It should be mentioned that
Foley, Ho\`{a}ng and Merkel's conjecture is stronger than Stanley's $\mathbf{(3+1)}$-free conjecture, since any claw-free incomparability graph is (claw, net)-free in view of the fact that the net is not an incomparability graph.

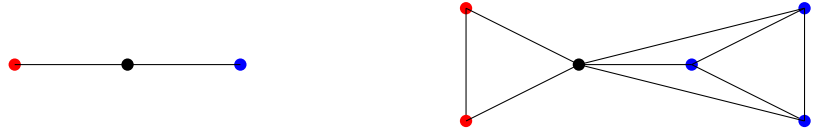
\begin{figure}[ht]
\centering
\begin{tikzpicture}[scale = 1.5]
    \fill[red] (-4,-0.5) circle (0.3ex);
    \fill (-3,-0.5) circle (0.3ex);
    \fill[blue] (-2,-0.5) circle (0.3ex);
    \draw (-4,-0.5) -- (-3,-0.5);
    \draw (-3,-0.5) -- (-2,-0.5);
    \fill[red] (0,0) circle (0.3ex);
    \draw (0,0) -- (0,-1);
    \fill[red] (0,-1) circle (0.3ex);
    \draw (0,0) -- (1,-0.5);
    \fill[blue] (2,-0.5) circle (0.3ex);
    \draw (0,-1) -- (1,-0.5);
    \fill[blue] (3,0) circle (0.3ex);
    \fill (1,-0.5) circle (0.3ex);
    \draw (1,-0.5) -- (2,-0.5);
    \draw (3,0) -- (1,-0.5);
    \draw (3,-1) -- (2,-0.5);
    \fill[blue] (3,-1) circle (0.3ex);
    \draw (2,-0.5) -- (3,0);
    \draw (1,-0.5) -- (3,-1);
    \draw (3,0) -- (3,-1);
\end{tikzpicture}
    \caption{The path $P_3$ and an expansion of $P_3$.}\label{ep3}
\end{figure}

In particular, Foley, Ho\`{a}ng and Merkel \cite{FHM19} studied the $e$-positivity of (claw, bull)-free graphs, which form a subfamily of (claw, net)-free graphs (the bull graph is depicted in Figure \ref{clawnetbull}). Hermosillo de la Maza, Jing and Masjoody \cite{HJM18} showed that each connected (claw, bull)-free graph is either a co-triangle-free graph, or an expansion of a path of length at least four or an expansion of a cycle of length at least six.
The co-triangle is the graph consisting of three isolated vertices and all co-triangle-free graphs are known to be $e$-positive \cite[Exercise 7.47]{StaEC2}. It remains to study the $e$-positivity of every expansion of a path or cycle. Recall that a graph is an \textit{expansion} of a path (resp.\ a cycle) if it is obtained from a path (resp.\ a cycle) by replacing every vertex $v_i$ by a complete graph $K_{n_i}$ ($n_i \ge 1$) while preserving the previous adjacency. For instance, Figure \ref{ep3} shows the path~$P_3$ and the expansion graph by replacing the three vertices with $K_2$, $K_1$, $K_3$ from left to right respectively. These expansion graphs are also called \textit{clan} graphs, see \cite{Rea81}.

Given a graph $G$ and a vertex $v$, repeatedly replacing $v$ by a $K_2$ for $n$ steps is equivalent to replacing $v$ by a $K_{n+1}$. This motivates Foley, Ho\`{a}ng and Merkel to focus on the expansion graph obtained by replacing $v$ by a $K_2$.
Applying such an operation is equivalent to adding a twin vertex $v'$ to $G$ such that $v,v'$ are adjacent and any other vertex of $G$ is adjacent to both of them or neither of them. For this reason we call it the twinning operation. We denote by $G_v$ the resulting graph. Foley, Ho\`{a}ng and Merkel proposed the following conjecture.

\begin{conj}\cite[Conjecture 23]{FHM19}\label{mconj}
If $G$ is $e$-positive, then so is $G_v$ for any vertex $v \in V(G)$.
\end{conj}

It is well known that an $e$-positive symmetric function must be $s$-positive. Thus the validness of Conjecture \ref{mconj} would imply the following corollary.

\begin{conjc}\label{mconj-coro}
If $G$ is $e$-positive, then $G_v$ is $s$-positive for any vertex $v \in V(G)$.
\end{conjc}

If Conjecture \ref{mconj} is true, then the $e$-positivity of (claw, bull)-free graphs would be established since paths and cycles are $e$-positive \cite{Sta95}. Unfortunately, it is not true in general, and even Conjollary \ref{mconj-coro} is not true, as will be shown in Section \ref{cx-e-pos} and Section \ref{cx-s-pos}.

The remaining of this paper is organized as follows. In Section \ref{back} we review some definitions and results on symmetric functions, which will be used in subsequent sections. In Section \ref{cx-e-pos} we show that Conjecture \ref{mconj} fails for a subfamily of tadpole graphs, as depicted in Figure \ref{tadgra}. In Section \ref{cx-s-pos}, we show that certain twin graph $G_v$ is not even $s$-positive for some $e$-positive graph $G$.
Some further directions inspired by Conjecture \ref{mconj} are discussed in Section \ref{sect-last}.

\section{Preliminaries}\label{back}

This section is devoted to giving an overview of some basic definitions and useful results on symmetric functions in commuting variables as well as those in non-commuting variables. Some results on chromatic symmetric functions in non-commuting variables are also covered here. The results given here play an important role in determining the $e$-positivity or $s$-positivity of the graphs to be studied in Sections \ref{cx-e-pos} and \ref{cx-s-pos}.

For notations and terminology on symmetric functions in commuting variables we shall follow Stanley \cite[Chapter 7]{StaEC2}.
Let ${\bf{x}}= \{x_1,x_2,\ldots\}$ be a set of countably infinite commuting indeterminates and $\mathbb{Q}[[\mathbf{x}]]$ be the algebra of formal power series over the rational field $\mathbb{Q}$ in $\mathbf{x}$. Then the \textit{algebra of symmetric functions} $\Lambda_{\mathbb{Q}}(\mathbf{x})$ is defined to be the subalgebra of $\mathbb{Q}[[\mathbf{x}]]$ consisting of all formal power series $f$ which are of bounded degree and satisfy
\[
	f(x_1,x_2,\ldots) = f(x_{\omega(1)},x_{\omega(2)},\ldots)
\]
for any permutation $\omega$ of positive integers.

The bases of $\Lambda_{\mathbb{Q}}(\textbf{x})$ are indexed by integer partitions. An \textit{integer partition} $\lambda$ of $n$ is a sequence of weakly decreasing positive integers $(\lambda_1,\lambda_2,\ldots,\lambda_{\ell})$ such that $\sum \lambda_i = n$ (if~$\lambda$ has only one part, we shall omit the parenthesis for convenience). Here $\ell = \ell(\lambda)$ is called the \textit{length} of $\lambda$. If $\lambda$ has $r_i$ parts equal to $i$, then we write $\lambda =\langle 1^{r_1}, 2^{r_2}, \ldots, n^{r_n}\rangle$ or $\lambda = (n^{r_n},(n-1)^{r_{n-1}},\ldots,1^{r_1})$, where terms with $r_i=0$ and the superscript $r_i=1$ may be omitted. We also use the notation $\lambda\vdash n$ to represent that $\lambda$ is a partition of $n$. For example, we may write $\lambda = (4,3,3,2,2,2) = \langle 2^3, 3^2, 4 \rangle = (4,3^2,2^3) \vdash 16$.

Three bases of $\Lambda_{\mathbb{Q}}(\textbf{x})$ are involved in this paper: the monomial symmetric functions~$m_{\lambda}$, the elementary symmetric functions $e_{\lambda}$, and the Schur functions $s_{\lambda}$. For any symmetric function $f\in\Lambda_{\mathbb{Q}}(\textbf{x})$ and a basis $\{u_{\lambda}\}$, we use the notation $[u_\lambda]f$ to denote the coefficient $c_\lambda$ in the expansion $f=\sum_{\lambda}c_\lambda u_{\lambda}$. In order to define $m_{\lambda}$, we first identify the sequence $(\lambda_1,\lambda_2,\ldots,\lambda_{\ell})$ with the infinite sequence $(\lambda_1,\lambda_2,\ldots,\lambda_{\ell},0,0,\ldots)$. Then
\[
	m_{\lambda} = \sum_{\alpha} x^{\alpha},
\]
where $\alpha = (\alpha_1,\alpha_2,\ldots)$ ranges over all distinct permutations of $\lambda$ and $x^{\alpha}$ represents the product $x_1^{\alpha_1}x_2^{\alpha_2}\cdots$. The elementary symmetric functions are defined by
\[
	e_{0} = 1, \quad e_{k} = \sum_{i_1 < i_2 < \cdots < i_k} x_{i_1}x_{i_2}\cdots x_{i_k} \text{ for } k \ge 1 \quad \text{and} \quad e_{\lambda} = e_{\lambda_1}e_{\lambda_2}\cdots e_{\lambda_{\ell}}.
\]

\begin{figure}
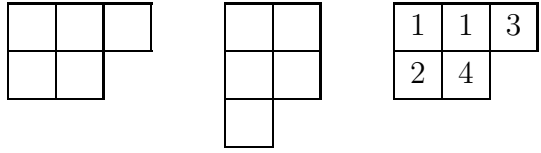

\begin{center}
\begin{ytableau}[]
& &\\
&
\end{ytableau} \qquad
\begin{ytableau}[]
 & \\
 & \\
 & \none
\end{ytableau} \qquad
\begin{ytableau}[]
1 & 1 & 3\\
2 & 4
\end{ytableau}
\end{center}
\caption{Examples for Young diagrams and Young tableaux.}\label{Young}
\end{figure}

A combinatorial definition of Schur functions relies on Young tableaux. For any partition $\lambda$, we can associate with it a left-justified diagram whose $i$-th row has $\lambda_i$ boxes, which is called the \textit{Young diagram} of $\lambda$. A \textit{filling} of a Young diagram is a way of putting numbers into the boxes. A \textit{semi-standard Young tableau} of shape $\lambda$ is defined to be a filling of the Young diagram of $\lambda$ such that the numbers in each row is weakly increasing from left to right and the numbers in each column is strictly increasing from top to bottom. If there are $c_1$ 1's, $c_2$ 2's, $\ldots$ in a Young tableau, then we say the \textit{type} of this tableau is $(c_1,c_2,\ldots)$. The \textit{conjugate partition} of $\lambda$, denoted by $\lambda'$, is the partition whose diagram is the transpose of the diagram of $\lambda$. For instance, the Young diagram of $\lambda = (3,2)$, that is $\lambda' = (2,2,1)$, and a semi-standard Young tableau of shape $\lambda$ and type $(2,1,1,1)$ are presented in Figure \ref{Young} from left to right. 
The Schur function $s_{\lambda}$ is defined by
\[
	s_{\lambda} = \sum_T x^T,
\]
where $T$ ranges over all semi-standard Young tableaux of shape $\lambda$, and $x^T = x_1^{c_1}x_2^{c_2}\cdots$ if $T$ is of type $(c_1,c_2,\ldots)$. The aforementioned three bases are transited via the Kostka numbers $K_{\lambda\mu}$, which is the number of semi-standard Young tableaux of shape $\lambda$ and type~$\mu$.

\begin{prop}\cite[p.311, p.335]{StaEC2}\label{smes}
We have
\begin{align}
	s_{\lambda} = \sum_{\mu} K_{\lambda,\mu} m_{\mu} \qquad \mbox{ and }\qquad 	e_{\mu} = \sum_{\lambda} K_{\lambda',\mu} s_{\lambda}.
\end{align}
\end{prop}

We would like to point out that the transition matrices $(K_{\lambda,\mu})_{\lambda,\mu\vdash n}$ could be written as a triangular matrix by suitably arranging the partitions of $n$ due to the following property satisfied by the
Kostka numbers.

\begin{prop}\cite[Proposition 7.10.5]{StaEC2}\label{Kosnum}
If $K_{\lambda,\mu} \neq 0$, then $\mu \le \lambda$ under dominance order, i.e., $\sum_{i=1}^{\ell(\mu)} \mu_i = \sum_{j=1}^{\ell(\lambda)} \lambda_j$ and $\sum_{i=1}^k \mu_i \le \sum_{i=1}^k \lambda_i$ for all $k \ge 1$. Moreover, $K_{\lambda,\lambda} = 1$.
\end{prop}

In general, there are no explicit formulas to compute the Kostka numbers $K_{\lambda\mu}$. But for $\mu=\langle 1^n\rangle$ and $\lambda\vdash n$, in which case $K_{\lambda\mu}$ is usually denoted by $f^{\lambda}$, there is the \textit{hook length formula}. In this case, the Young tableaux are called \textit{standard} Young tableaux. Recall that the \textit{hook} of a box $b$ in a Young diagram consists of $b$, the boxes to the right of $b$ in the same row, and the boxes below $b$ in the same column. The \textit{hook length} of $b$, denoted by $h(b)$, is the number of boxes in its hook. The hook length formula is stated as follows.

\begin{prop}\cite{FRT54}\label{hlf}
If $\lambda\vdash n$, then
\[
	f^{\lambda} = \frac{n!}{\prod_{b} h(b)},
\]
where $b$ ranges over all boxes of the Young diagram of $\lambda$.
\end{prop}

Next, we review some results on chromatic symmetric functions from \cite{Sta95}. The following proposition gives the relationship between the chromatic symmetric function of a graph $G$ and its connected components.

\begin{prop}\cite[Proposition 2.3]{Sta95}\label{prop-un}
Let $G \uplus H$ denote the disjoint union of graphs $G$ and $H$. Then
\[X_{G \uplus H}=X_GX_H.\]
\end{prop}

In Section \ref{cx-s-pos} we will use the monomial expansion of chromatic symmetric functions, for which Stanley gave a combinatorial interpretation of the coefficients by stable partitions. Recall that a \textit{stable partition} of a graph $G$ is a set partition of its vertex set $V(G)$ such that each block is a stable set, that is, any two vertices in the same block are not adjacent. The \textit{type} of a stable partition is defined to be the integer partition whose parts are the block sizes. Then we have the following result.

\begin{prop}\cite[Proposition 2.4]{Sta95}\label{stanley-lemma}
Let $G$ be a graph of $d$ vertices and $a_{\lambda}$ be the number of stable partitions of $G$ of type $\lambda$. Then
\[
    X_G = \sum_{\lambda} a_{\lambda}\tilde{m}_{\lambda},
\]
where $\tilde{m}_{\lambda} = r_1!r_2!\cdots r_d!m_{\lambda}$ if $\lambda =\langle 1^{r_1},2^{r_2},\ldots, d^{r_d}\rangle$.
\end{prop}

To study Stanley's chromatic symmetric functions, Gebhard and Sagan \cite{GS01} made a systematic study of the algebra of symmetric functions in non-commuting variables, for which we continue to use ${\bf{x}}= \{x_1,x_2,\ldots\}$ to represent the indeterminants but we require that $x_ix_j\neq x_jx_i$ for any $i\neq j$. Let $G$ be a graph with vertices labeled $v_1,v_2,\ldots,v_d$ in a fixed order. Gebhard and Sagan \cite{GS01}  defined the \textit{chromatic symmetric function in non-commuting variables} of $G$ as
\[Y_G=\sum_{\kappa}x_{\kappa(v_1)}x_{\kappa(v_2)}\cdots x_{\kappa(v_d)},\]
where the sum ranges over all proper colorings $\kappa$ of $G$.
Note that $Y_G$ depends on the labeling of $G$, and hence when using $Y_G$ we shall always deal with labeled graphs. Moreover, it is clear that $Y_G$ becomes $X_G$ if we let the variables commute.

In order to study the $e$-positivity of $X_G$, it is desirable to expand $Y_G$
in terms of  the elementary symmetric functions in non-commuting variables. Different from the classical case, such
elementary symmetric functions are indexed by set partitions. Let $\Pi_d$ denote the collection of set partitions $\pi$ of the set $[d]:=\{1,2,\ldots,d\}$. Given $\pi\in\Pi_d$, the \textit{elementary symmetric function in non-commuting variables} $e_\pi$ is defined by
\[e_\pi=\sum_{(i_1,i_2,\ldots,i_d)} x_{i_1}x_{i_2}\cdots x_{i_d},\]
where the sum ranges over all sequences $(i_1,i_2,\ldots,i_d)$ of positive integers such that $i_j\neq i_k$ if $j$ and $k$ are in the same block of $\pi$. For $\pi\in\Pi_d$ we define the \textit{type} of $\pi$, denoted $\lambda(\pi)$, to be the integer partition of $d$ whose parts are the sizes of the blocks of $\pi$. Clearly, if the variables are allowed to commute, then $e_{\pi}$ becomes $\lambda(\pi)! e_{\lambda(\pi)}$, where $\lambda(\pi)! = \prod_{i} \lambda(\pi)_i !$. Therefore, given $Y_G = \sum_{\pi} c_{\pi} e_{\pi}$, if these set partitions could be divided into disjoint sets $D_1,D_2,\ldots$ such that all the set partitions in the same $D_i$ are of the same type and $\sum_{\pi \in D_i} c_{\pi} \ge 0$, then
\begin{equation*}\label{xd}
X_G = \sum_{D_i} \left(\sum_{\pi \in D_i} c_{\pi}\right) \lambda(\pi)! e_{\lambda(\pi)}
\end{equation*}
is $e$-positive. In light of this, Gebhard and Sagan introduced the following equivalence relation on set partitions.
Given two set partitions $\sigma$ and $\tau$ of $[d]$ and an integer $i\in [d]$, let $B_{\sigma,i}$ and $B_{\tau,i}$ denote the block of $\sigma$ containing $i$ and the block of $\tau$ containing $i$, respectively. We say that $\sigma$ is congruent to $\tau$ modulo $i$, denoted by $\sigma\equiv_i\tau$, if
\[\lambda(\sigma)=\lambda(\tau)\ \mbox{and}\ |B_{\sigma,i}|=|B_{\tau,i}|.\]
It is easy to verify that the relation $\equiv_i$ is an equivalence relation on $\Pi_d$.
Further, we can extend this equivalence relation to the set of elementary symmetric functions in non-commuting variables by letting
\[e_\sigma\equiv_i e_\tau\ \iff\ \sigma\equiv_i\tau.\]
Let $(\tau)$ and $e_{(\tau)}$ denote the equivalence classes of $\tau$ and $e_\tau$, respectively. If
$$Y_G = \sum_{\pi\in \Pi_d} c_{\pi} e_{\pi},$$
then we can collect the equivalent terms and rewrite
\begin{align}\label{mod-expansion}
Y_G \equiv_i \sum_{(\tau)}c_{(\tau)} e_{(\tau)},
\end{align}
where the sum ranges over all equivalence classes of $\Pi_d$ and $c_{(\tau)}=\sum_{\sigma\in(\tau)} c_\sigma$.
If all $c_{(\tau)}$ in \eqref{mod-expansion} are non-negative, then we say that $G$ is ($e$)-\textit{positive} modulo $i$. It is easy to see that the $(e)$-positivity of $Y_G$ would imply the $e$-positivity of $X_G$.

In Section \ref{cx-e-pos} we will establish the $(e)$-positivity of tadpole graphs based on two results given by Gebhard and Sagan \cite{GS01}. The first one is on the $(e)$-positivity of a cycle, whose $e$-positivity was already proved by Stanley \cite{Sta95}.

\begin{prop}\cite[Proposition 6.8]{GS01}\label{C-thm}
For all $d\geq 2$, if a cycle $C_d$ is labeled as $V(C_d) = \{v_1,\ldots,v_d\}$, $E(C_d) = \{v_1 v_2,v_2 v_3,\ldots,v_{d-1} v_d,v_d v_1\}$, then $Y_{C_d}$ is ($e$)-positive modulo $d$.
\end{prop}

The second result allows us to construct new $(e)$-positive graphs from old ones.

\begin{prop}\cite[Lemma 7.5]{GS01}\label{G+K-thm}
Given $n,d\geq 1$, let $G$ be a graph with vertex set $\{v_1,v_2,\ldots,v_d\}$, and let $G+K_n$ be the graph with vertex set $V(G+K_n)=V(G)\cup\{v_{d+1},\ldots,v_{d+n-1}\}$ and edge set $E(G+K_n)=E(G)\cup\{e=v_iv_j:i,j\in [d,d+n-1]\}$.
If $Y_G$ is ($e$)-positive modulo $d$, then $Y_{G+K_n}$ is ($e$)-positive modulo $d+n-1$.
\end{prop}

\section{Counterexamples to Conjecture \ref{mconj}} \label{cx-e-pos}

The aim of this section is to construct counterexamples to Conjecture \ref{mconj} from tadpole graphs. Recall that the $(a,b)$-\textit{tadpole graph} $T^{\langle a,b\rangle}$ ($a \ge 3, b \ge 1$) is obtained by connecting an $a$-vertex cycle $C_{a}$ and a $b$-vertex path $P_{b}$ with an edge, as shown in Figure \ref{tadgra}. In this section we first prove the $e$-positivity of $T^{\langle a,b\rangle}$ for any $a \ge 3, b \ge 1$. Next we show that, for any $b\geq 1$, the twin graph $T^{\langle 4,b\rangle}_{v_4}$ is not $e$-positive, where $v_{4}$ is the unique vertex of degree 3 in $T^{\langle 4,b\rangle}$. This implies that Conjecture \ref{mconj} fails for $\{T^{\langle 4,b \rangle}\, |\, b \ge 1\}$.

\begin{figure}[ht]
\centering
\begin{tikzpicture}
[place/.style={thick,fill=black!100,circle,inner sep=0pt,minimum size=1mm,draw=black!100},rotate=180]
\node [place, label=above:$v_{a+b-1}$] (v3) at (-1,7) {};
\node [place, label=above:$v_{a+b}$] (v2) at (-2,7) {};
\draw [thick] (v2) -- (v3);
\node [place, label=above:$v_{a+1}$] (v4) at (0,7) {};
\node [place, label=left:$v_{a}$] (v5) at (1,7) {};
\draw [thick] (v4) -- (v5);
\node at (-0.5,7) {$\cdots$};
\draw [thick] (2,6) arc (-90:-270:10mm);
\draw [thick] (2,6) [densely dashed] arc (-90:90:10mm);
\node [place, label=above:$v_{1}$] at (2,6) {};
\node [place, label=below:$v_{a-1}$] at (2,8) {};
\end{tikzpicture}
\caption{The tadpole graph $T^{\langle a,b \rangle}.$}\label{tadgra}
\end{figure}

Let us first prove the $e$-positivity of tadpole graphs by using Propositions \ref{C-thm} and \ref{G+K-thm}.

\begin{thm}\label{tadepos}
For any $a\geq 3$ and $b\geq 1$, the tadpole graph $T^{\langle a,b \rangle}$ is $e$-positive.
\end{thm}

\pf We use the labeling of $T^{\langle a,b \rangle}$ in Figure \ref{tadgra}. Note that this labeling of $C_{a}$ coincides with that of Proposition \ref{C-thm}, which implies that $C_{a}$ is ($e$)-positive modulo $a$. Next, we join a $K_2$, which is labeled as $v_{a}v_{a+1}$, to the vertex $v_{a}$. Then by Proposition \ref{G+K-thm} this new graph is $(e)$-positive modulo $a+1$. Now we join $v_{a+1}v_{a+2}$ to $v_{a+1}$ and the resulting graph is $(e)$-positive modulo $a+2$. Continuing this process, we will finally get that $T^{\langle a,b \rangle}$ is $(e)$-positive modulo $a+b$, and hence $T^{\langle a,b \rangle}$ is $e$-positive. \qed

In the remaining part of this section, we concentrate on proving that $T^{\langle 4,b \rangle}_{v_4}$ is not $e$-positive for any $b \ge 1$.
To this end, we will expand $X_{T^{\langle 4,b \rangle}_{v_4}}$ in terms of the chromatic symmetric functions of paths.
For convenience of illustration, we shall keep using the labeling of vertices of $T^{\langle 4,b \rangle}_{v_4}$ in Figure \ref{tadgra'}, where the twin vertices are labeled as $v_4$ and $v'_4$.

\begin{figure}[ht]
\centering
\begin{tikzpicture}
[place/.style={thick,fill=black!100,circle,inner sep=0pt,minimum size=1mm,draw=black!100}, sample/.style={thick,fill=blue!100,circle,inner sep=0pt,minimum size=1mm,draw=blue!100},rotate=180]
\node [place, label=above:$v_{b+4}$] (v1) at (-2,4.5) {};
\node [place, label=above:$v_{b+3}$] (v6) at (-1,4.5) {};
\node at (-0.5,4.5) {$\cdots$};
\node [place, label=above:$v_5$] (v7) at (0,4.5) {};
\node [sample, label=below:$v_{4}$] (v8) at (1,5) {};
\node [place, label=below:$v_{3}$] (v10) at (2,5.5) {};
\node [place, label=above:$v_{1}$] (v11) at (2,3.5) {};
\draw [thick]  (v1) edge (v6);
\draw [thick]  (2,3.5) arc (-90:90:10mm);
\node [sample, label=above:$v_{4}'$] (v9) at (1,4) {};
\draw [thick] (v7) -- (v8) -- (v9) -- (v7);
\draw [thick] (v10) -- (v8) -- (v11);
\draw [thick] (v10) -- (v9) -- (v11);
\node [place, label=left:$v_{2}$] at (3,4.5) {};
\node [blue] at (0.5,5) {$e_2$};
\node [blue] at (0.5,4) {$e_1$};
\node [blue] at (0.75,4.5) {$e_3$};
\end{tikzpicture}
\caption{The graph $T^{\langle 4,b \rangle}_{v_{4}}.$}\label{tadgra'}
\end{figure}

For our desired expansion of $X_{T^{\langle 4,b \rangle}_{v_4}}$, we shall use two results due to  Orellana and Scott \cite{OS14}, which enable one to express the chromatic symmetric function of a graph by that of other ``simpler'' graphs. Given $S \subseteq E(G)$ and $T \subseteq \{v_iv_j \mid v_iv_j \not\in E(G)\}$, let $(G - S)\cup T$ denote the graph obtained from $G$ by removing the edges in $S$ and adding the edges in $T$. Orellana and Scott's results are stated as follows.

\begin{prop}\cite[Theorem 3.1]{OS14}\label{dec-thm}
If $G=(V,E)$ is a graph with $e_1, e_2, e_3\in E$ forming a triangle, then
\[X_G=X_{G-\{e_1\}}+X_{G-\{e_2\}}-X_{G-\{e_1,e_2\}}.\]
\end{prop}

\begin{prop}\cite[Corollary 3.2]{OS14}\label{dec-cor}
Let $G=(V,E)$ be a graph and let $v_1,v_2,v_3\in V$ satisfying $e_1=v_1v_3\in E$, $e_2=v_2v_3\in E$ and $e_3=v_1v_2\not\in E$.
Then
\[X_G=X_{(G-\{e_1\})\cup\{e_3\}}+X_{G-\{e_2\}}-X_{(G-\{e_1,e_2\})\cup\{e_3\}}.\]
\end{prop}

Now we are in the position to state the expansion of $X_{T^{\langle 4,b\rangle}_{v_{4}}}$.

\begin{prop}\label{PC-lem}
For any $b \geq 1$ we have
\begin{align}\label{t'}
X_{T^{\langle 4,b \rangle}_{v_4}}=&20 X_{P_{b+5}}+2 e_1 X_{P_{b+4}}-16 e_2 X_{P_{b+3}}-(2 e_{(2,1)} +42 e_3 ) X_{P_{b+2}} \notag \\
&-(56e_4+4e_{(2,2)}+4e_{(3,1)}) X_{P_{b+1}}  -(6e_{(4,1)}+4e_{(3,2)}+50e_5 ) X_{P_{b}}.
\end{align}
\end{prop}

\begin{figure}[ht]
\centering
\begin{tikzpicture}
[place/.style={thick,fill=black!100,circle,inner sep=0pt,minimum size=1mm,draw=black!100}, sample/.style={thick,fill=blue!100,circle,inner sep=0pt,minimum size=1mm,draw=blue!100},rotate=180]
\node [place, label=above:$v_{b+4}$] (v1) at (-2,4.5) {};
\node [place, label=above:$v_{b+3}$] (v6) at (-1,4.5) {};
\node at (-0.5,4.5) {$\cdots$};
\node [place, label=above:$v_5$] (v7) at (0,4.5) {};
\node [sample, label=below:$v_{4}$] (v8) at (1,5) {};
\node [place, label=below:$v_{3}$] (v10) at (2,5.5) {};
\node [place, label=above:$v_{1}$] (v11) at (2,3.5) {};
\draw [thick]  (v1) edge (v6);
\draw [thick]  (2,3.5) arc (-90:90:10mm);
\node [sample, label=above:$v_{4}'$] (v9) at (1,4) {};
\draw [thick] (v7) -- (v8) -- (v9);
\draw [thick] (v10) -- (v8) -- (v11);
\draw [thick] (v10) -- (v9) -- (v11);
\node [place, label=left:$v_{2}$] at (3,4.5) {};
%\node [blue] at (0.5,5) {$e_2$};
%\node [blue] at (0.75,4.5) {$e_3$};

\node [place] (v2) at (-6.5,5) {};
\node [place] (v3) at (-6.5,4) {};
\node [place] (v4) at (-5.5,3.5) {};
\node [place] (v1) at (-5.5,5.5) {};
\node [place] at (-4.5,4.5) {};
\draw [thick] (-5.5,3.5) arc (-90:90:1);
\draw [thick] (v4) -- (v3) -- (v2) -- (v1);
\draw [thick] (v1) -- (v3);
\draw [thick] (v2) -- (v4);
\node at (0,6.5) {$H_b$};
\node at (-5.5,6.5) {$S$};
\end{tikzpicture}
\caption{The graphs $H_b$ and $S$.}\label{HandS}
\end{figure}
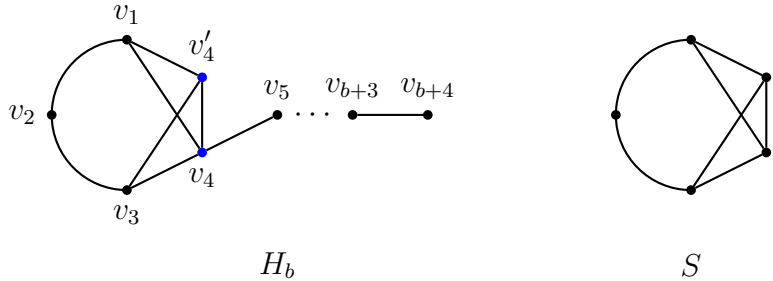

\pf In Proposition \ref{dec-thm} we take $G=T^{\langle 4,b \rangle}_{v_4}$, $e_1=v_4'v_5$, $e_2=v_4v_5$ and $e_3=v_4v_4'$, as labeled in Figure \ref{tadgra'}. It is easy to see that $G-\{e_1\}$, denoted by $H_b$, and $G-\{e_2\}$ are isomorphic, and that $G-\{e_1,e_2\}$ is the disjoint union of a $b$-vertex path $P_b$ and the ``cycle part'' $S$ of $T^{\langle 4,b \rangle}_{v_{4}}$, as shown in Figure \ref{HandS}. Thus by Proposition \ref{dec-thm} and Proposition \ref{prop-un},
\begin{align}\label{xt1}
X_{T^{\langle 4,b\rangle}_{v_{4}}}=2X_{H_b}-X_{P_b}X_{S}.
\end{align}

By direct computation we get
\begin{align}\label{xs}
	X_S = 4e_{(3, 2)} + 6e_{(4, 1)} + 50e_{5},
\end{align}
and hence we only need to expand $X_{H_b}$.

To decompose $X_{H_b}$, we label $v_3v_4$, $v_3v_4'$, $v_4v_4'$ by $e_1$, $e_2$ and $e_3$, respectively, as shown in Figure \ref{gra-h1}. By applying Proposition \ref{dec-thm} we deduce that
\begin{align}\label{equ-H-1}
X_{H_b}=X_{H_b-\{v_3v_4\}}+X_{H_b-\{v_3v_4'\}}-X_{H_b-\{v_3v_4,v_3v_4'\}}.
\end{align}
%where $H'_{2,3}$, $H'_{1,3}$ and $H'_3$ are obtained by deleting $e_1$, $e_2$ and $\{e_1,e_2\}$ in $H$, respectively.

\begin{figure}[ht]
\centering
\begin{tikzpicture}
[place/.style={thick,fill=black!100,circle,inner sep=0pt,minimum size=1mm,draw=black!100}, sample/.style={thick,fill=blue!100,circle,inner sep=0pt,minimum size=1mm,draw=blue!100},rotate=180]
\node [place, label=above:$v_{b+4}$] (v1) at (-2,4.5) {};
\node [place, label=above:$v_{b+3}$] (v6) at (-1,4.5) {};
\node at (-0.5,4.5) {$\cdots$};
\node [place, label=above:$v_5$] (v7) at (0,4.5) {};
\node [sample, label=below:$v_{4}$] (v8) at (1,5) {};
\node [place, label=below:$v_{3}$] (v10) at (2,5.5) {};
\node [place, label=above:$v_1$] (v11) at (2,3.5) {};
\draw [thick]  (v1) edge (v6);
\draw [thick]  (2,3.5) arc (-90:90:10mm);
\node [sample, label=above:$v_{4}'$] (v9) at (1,4) {};
\draw [thick] (v7) -- (v8) -- (v9);
\draw [thick] (v10) -- (v8) -- (v11);
\draw [thick] (v10) -- (v9) -- (v11);
\node [place, label=left:$v_{2}$] at (3,4.5) {};
\node [blue] at (1.5,5.5) {$e_1$};
\node [blue] at (1.75,4.75) {$e_2$};
\node [blue] at (0.75,4.5) {$e_3$};
\end{tikzpicture}
\caption{
%An illustration of Lemma \ref{PC'-2}
The graph $H_b$.}\label{gra-h1}
\end{figure}

\newpage

Next we calculate $X_{H_b-\{v_3v_4\}}$, $X_{H_b-\{v_3v_4'\}}$ and $X_{H_b-\{v_3v_4,v_3v_4'\}}$ individually. We shall first expand each of them as a linear combination of some ``intermediate'' graphs, and then give the final expansion in terms of $X_{P_n}$.

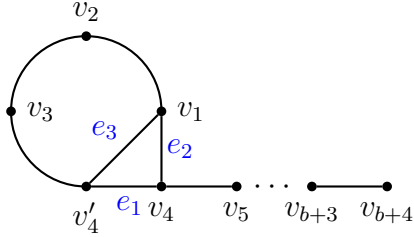
\begin{figure}[h]
\centering
\begin{tikzpicture}
[place/.style={thick,fill=black!100,circle,inner sep=0pt,minimum size=1mm,draw=black!100}, sample/.style={thick,fill=blue!100,circle,inner sep=0pt,minimum size=1mm,draw=blue!100},rotate=180]
\node [place, label=below:$v_{b+4}$] (v1) at (-2,4.5) {};
\node [place, label=below:$v_{b+3}$] (v6) at (-1,4.5) {};
\node at (-0.5,4.5) {$\cdots$};
\node [place, label=below:$v_5$] (v7) at (0,4.5) {};
\node [place, label=right:$v_{1}$] (v8) at (1,3.5) {};
\node [place, label=above:$v_{2}$] (v10) at (2,2.5) {};
\node [place, label=below:$v_{4}'$] (v11) at (2,4.5) {};
\draw [thick]  (v1) edge (v6);
\draw [thick]  (2,4.5) arc (90:-180:10mm);
\node [place, label=below:$v_{4}$] (v9) at (1,4.5) {};
\draw [thick] (v8) -- (v9) -- (v7);
\draw [thick] (v8) -- (v11);
\draw [thick]  (v9) -- (v11);
\node [place, label=right:$v_{3}$] at (3,3.5) {};
\node [blue] at (1.425,4.75) {$e_1$};
\node [blue] at (0.75,4) {$e_2$};
\node [blue] at (1.75,3.75) {$e_3$};
\end{tikzpicture}
\caption{
%An illustration of Lemma \ref{PC'-2}
$H_b-\{v_3v_4\}$.}\label{gra-h2}
\end{figure}

For $X_{H_b-\{v_3v_4\}}$, we label $v_4v_4'$, $v_1v_4$, $v_1v_4'$ by $e_1$, $e_2$ and $e_3$ respectively, as shown in Figure \ref{gra-h2}. Then by Proposition \ref{dec-thm} we have
\begin{align}\label{equ-H-2}
X_{H_b-\{v_3v_4\}}=2X_{T^{\langle 4,b+1 \rangle}}-X_{P_{b+1}}X_{C_4}.
\end{align}

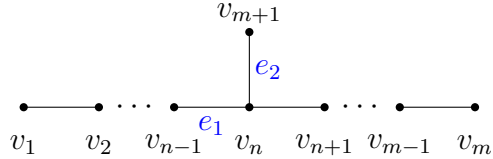
\begin{figure}[ht]
\centering
\begin{tikzpicture}
%[scale = 1.5]
    \fill (0,0) circle (0.3ex);
    \draw (0,0) -- (1,0);
    \fill (1,0) circle (0.3ex);
    \node at (1.5,0) {$\cdots$};
    \draw (2,0) -- (3,0);
    \fill (2,0) circle (0.3ex);
    \fill (3,0) circle (0.3ex);
    \node at (4.5,0) {$\cdots$};
    \fill (3,1) circle (0.3ex);
    \draw (3,0) -- (3,1);
    \fill (4,0) circle (0.3ex);
    \draw (3,0) -- (4,0);
    \fill (5,0) circle (0.3ex);
    \draw (5,0) -- (6,0);
    \fill (6,0) circle (0.3ex);
    \node at (0,-0.5) {$v_1$};
    \node at (1,-0.5) {$v_2$};
    \node at (2,-0.5) {$v_{n-1}$};
    \node at (3,-0.5) {$v_n$};
    \node at (4,-0.5) {$v_{n+1}$};
    \node at (5,-0.5) {$v_{m-1}$};
    \node at (6,-0.5) {$v_{m}$};
    \node [blue] at (2.5,-0.25) {$e_1$};
    \node [blue] at (3.25,0.5) {$e_2$};
    \node at (3,1.25) {$v_{m+1}$};
\end{tikzpicture}
    \caption{The graph $P_{m,n}$.}\label{gra-pmn}
\end{figure}

We next compute $X_{H_b-\{v_3v_4'\}}$. For $n\leq m$, we define $P_{m,n}$ to be the graph obtained by adding a leaf vertex to the $n$-th vertex of $P_m$, as shown in Figure \ref{gra-pmn}.

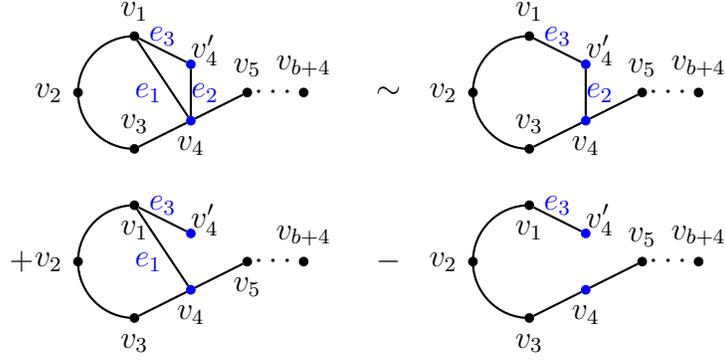
\begin{figure}[ht]
  \centering
\begin{tikzpicture}
[place/.style={thick,fill=black!100,circle,inner sep=0pt,minimum size=1mm,draw=black!100}, sample/.style={thick,fill=blue!100,circle,inner sep=0pt,minimum size=1mm,draw=blue!100},scale=0.75,rotate=180]
\node [place, label=above:$v_{b+4}$] (v6) at (-1,4.5) {};
\node at (-0.5,4.5) {$\cdots$};
\node [place, label=above:$v_5$] (v7) at (0,4.5) {};
\node [sample, label=below:$v_{4}$] (v8) at (1,5) {};
\node [place, label=below:$v_{3}$] (v10) at (2,5.5) {};
\node [place, label=below:$v_{1}$] (v11) at (2,3.5) {};
\draw [thick]  (2,3.5) arc (-90:90:10mm);
\node [sample] (v9) at (1,4) {};
\node at (0.75,3.75) {$v_{4}'$};
\draw [thick] (v7) -- (v8);
\draw [thick] (v10) -- (v8);
\draw [thick] (v9) -- (v11);
\node [place, label=left:$v_{2}$] at (3,4.5) {};
\node [blue] at (1.5,3.5) {$e_3$};

\node at (4.5,1.5) {$\sim$};
\node [place, label=above:$v_{b+4}$] at (6,4.5) {};
\node at (6.5,4.5) {$\cdots$};
\node [place, label=below:$v_5$] (v1) at (7,4.5) {};
\node [sample, label=below:$v_{4}$] (v2) at (8,5) {};
\node [place, label=below:$v_{3}$] (v3) at (9,5.5) {};
\node [place, label=left:$v_{2}$] at (10,4.5) {};
\node [place, label=below:$v_{1}$] (v4) at (9,3.5) {};
\node [sample] (v5) at (8,4) {};
\node at (7.75,3.75) {$v_{4}'$};
\draw [thick]  (9,3.5) arc (-90:90:10mm);
\draw [thick] (v1) -- (v2) -- (v3);
\draw [thick] (v5) -- (v4) -- (v2);
\node [blue] at (8.75,4.5) {$e_1$};
\node [blue] at (8.5,3.5) {$e_3$};

\node at (11,4.5) {$+$};
\node [place, label=above:$v_{b+4}$]  at (-1,1.5) {};
\node at (-0.5,1.5) {$\cdots$};
\node [place, label=above:$v_5$] (v12) at (0,1.5) {};
\node [sample, label=below:$v_{4}$] (v13) at (1,2) {};
\node [place, label=above:$v_{3}$] (v16) at (2,2.5) {};
\node [place, label=left:$v_{2}$] at (3,1.5) {};
\node [place, label=above:$v_{1}$] (v15) at (2,0.5) {};
\draw [thick]  (2,0.5) arc (-90:90:10mm);
\node [sample] (v14) at (1,1) {};
\node at (0.75,0.75) {$v_{4}'$};
\draw [thick] (v13) -- (v16);
\draw [thick] (v12) -- (v13) -- (v14) -- (v15);
\node [blue]at (0.75,1.5) {$e_2$};
\node [blue]at (1.5,0.5) {$e_3$};

\node at (4.5,4.5) {$-$};
\node [place, label=above:$v_{b+4}$] at (6,1.5) {};
\node at (6.5,1.5) {$\cdots$};
\node [place, label=above:$v_5$] (v17) at (7,1.5) {};
\node [sample, label=below:$v_{4}$] (v18) at (8,2) {};
\node [place, label=above:$v_{3}$] (v19) at (9,2.5) {};
\node [place, label=left:$v_{2}$] at (10,1.5) {};
\node [place, label=above:$v_{1}$] (v20) at (9,0.5) {};
\draw [thick] (9,0.5) node {} arc (-90:90:10mm);
\node [sample] (v21) at (8,1) {};
\node at (7.75,0.75) {$v_{4}'$};
\draw [thick] (v17) -- (v18) -- (v19);
\draw [thick] (v21) edge (v20);
\node [blue] at (8.5,0.5) {$e_3$};
\draw [thick] (v21) -- (v18) -- (v20);
\node [blue] at (7.75,1.5) {$e_2$};
\node [blue] at (8.75,1.5) {$e_1$};
\end{tikzpicture}
\caption{Decomposition of $H_b-\{v_3v_4'\}$.}\label{gra-h3}
\end{figure}

Then by Proposition \ref{dec-thm} and Proposition \ref{dec-cor}, we could expand $X_{H_b-\{v_3v_4'\}}$ as
\begin{align}\label{equ-H-3}
X_{H_b-\{v_3v_4'\}} =& X_{T^{\langle 5,b \rangle}} + (X_{T^{\langle 5,b \rangle}} + X_{T^{\langle 4,b \rangle}} \cdot e_1 - X_{P_{b+4,4}}) - X_{P_{b+5}} \notag\\
=&2X_{T^{\langle 5,b \rangle}}+e_1X_{T^{\langle 4,b \rangle}}-X_{P_{b+4,4}}-X_{P_{b+5}}.
\end{align}
The procedure is illustrated in Figure \ref{gra-h3} and Figure \ref{gra-h4}, where for four graphs $G_1,G_2,G_3,G_4$, $G_1\sim G_2 + G_3 - G_4$ means $X_{G_1} = X_{G_2} + X_{G_3} - X_{G_4}$.

\newpage

\begin{figure}[ht]
\centering
\begin{tikzpicture}
[place/.style={thick,fill=black!100,circle,inner sep=0pt,minimum size=1mm,draw=black!100}, sample/.style={thick,fill=blue!100,circle,inner sep=0pt,minimum size=1mm,draw=blue!100},scale=0.75,rotate=180]
\node [place, label=above:$v_{b+4}$] (v6) at (-1,4.5) {};
\node at (-0.5,4.5) {$\cdots$};
\node [place, label=above:$v_5$] (v7) at (0,4.5) {};
\node [sample, label=below:$v_{4}$] (v8) at (1,5) {};
\node [place, label=below:$v_{3}$] (v10) at (2,5.5) {};
\node [place, label=below:$v_{1}$] (v11) at (2,3.5) {};
\draw [thick]  (2,3.5) arc (-90:90:10mm);
\node [sample] (v9) at (1,4) {};
\node at (0.75,3.75) {$v_{4}'$};
\draw [thick] (v7) -- (v8);
\draw [thick] (v10) -- (v8);
\draw [thick] (v8) edge (v9);
\node [place, label=left:$v_{2}$] at (3,4.5) {};
\node [blue] at (0.75,4.5) {$e_3$};

\node at (4.5,1.5) {$\sim$};
\node [place, label=above:$v_{b+4}$] at (6,4.5) {};
\node at (6.5,4.5) {$\cdots$};
\node [place, label=above:$v_5$] (v1) at (7,4.5) {};
\node [sample, label=below:$v_{4}$] (v2) at (8,5) {};
\node [place, label=below:$v_{3}$] (v3) at (9,5.5) {};
\node [place, label=left:$v_{2}$] at (10,4.5) {};
\node [place, label=below:$v_{1}$] (v4) at (9,3.5) {};
\node [sample] (v5) at (8,4) {};
\node at (7.75,3.75) {$v_{4}'$};
\draw [thick]  (9,3.5) arc (-90:90:10mm);
\draw [thick] (v1) -- (v2) -- (v3);
\node [blue] at (8.75,4.5) {$e_1$};
\draw [thick] (v4) edge (v2);

\node at (11,4.5) {$+$};
\node [place, label=above:$v_{b+4}$]  at (-1,1.5) {};
\node at (-0.5,1.5) {$\cdots$};
\node [place, label=above:$v_5$] (v12) at (0,1.5) {};
\node [sample, label=below:$v_{4}$] (v13) at (1,2) {};
\node [place, label=above:$v_{3}$] (v16) at (2,2.5) {};
\node [place, label=left:$v_{2}$] at (3,1.5) {};
\node [place, label=above:$v_{1}$] (v15) at (2,0.5) {};
\draw [thick]  (2,0.5) arc (-90:90:10mm);
\node [sample, label=above:$v_{4}'$] (v14) at (1,1) {};
\draw [thick] (v13) -- (v16);
\draw [thick] (v12) -- (v13);
\draw [thick] (v15) -- (v14) -- (v13);
\node [blue] at (1.5,0.5) {$e_2$};
\node [blue] at (0.75,1.5) {$e_3$};

\node at (4.5,4.5) {$-$};
\node [place, label=above:$v_{b+4}$] at (6,1.5) {};
\node at (6.5,1.5) {$\cdots$};
\node [place, label=above:$v_5$] (v17) at (7,1.5) {};
\node [sample, label=below:$v_{4}$] (v18) at (8,2) {};
\node [place, label=above:$v_{3}$] (v19) at (9,2.5) {};
\node [place, label=left:$v_{2}$] at (10,1.5) {};
\node [place, label=above:$v_{1}$] (v20) at (9,0.5) {};
\draw [thick] (9,0.5) arc (-90:90:10mm);
\node [sample, label=above:$v_{4}'$] (v21) at (8,1) {};
\draw [thick] (v17) -- (v18) -- (v19);
\node [blue] at (8.5,0.5) {$e_2$};
\draw [thick] (v21) -- (v20) -- (v18);
\node [blue] at (8.75,1.5) {$e_1$};

\end{tikzpicture}
\caption{Decomposition of $H_b-\{v_3v_4',\,v_4v_4'\}$, the third graph in Figure \ref{gra-h3}.}\label{gra-h4}
\end{figure}

\begin{figure}[ht]
\centering
\begin{tikzpicture}
[place/.style={thick,fill=black!100,circle,inner sep=0pt,minimum size=0.75mm,draw=black!100}, sample/.style={thick,fill=blue!100,circle,inner sep=0pt,minimum size=0.75mm,draw=blue!100},scale=0.75,rotate=180]
\node [place, label=above:$v_{b+4}$] (v6) at (-1,4.5) {};
\node at (-0.5,4.5) {$\cdots$};
\node [place, label=above:$v_5$] (v7) at (0,4.5) {};
\node [sample, label=below:$v_{4}$] (v8) at (1,5) {};
\node [place, label=below:$v_{3}$] (v10) at (2,5.5) {};
\node [place, label=above:$v_{1}$] (v11) at (2,3.5) {};
\draw [thick] (2,3.5) arc (-90:90:10mm);
\node [sample] (v9) at (1,4) {};
\draw [thick] (v7) -- (v8) -- (v11);
\node [place, label=left:$v_{2}$] at (3,4.5) {};
\node at (0.75,3.75) {$v_{4}'$};
\node [blue] at (1.75,4.5) {$e_3$};

\node at (4.5,4.5) {$-$};
\node [place, label=above:$v_{b+4}$] at (6,4.5) {};
\node at (6.5,4.5) {$\cdots$};
\node [place, label=above:$v_5$] (v1) at (7,4.5) {};
\node [sample, label=below:$v_{4}$] (v2) at (8,5) {};
\node [place, label=below:$v_{3}$] (v3) at (9,5.5) {};
\node [place, label=left:$v_{2}$] at (10,4.5) {};
\node [place, label=above:$v_{1}$] (v4) at (9,3.5) {};
\node [sample] (v5) at (8,4) {};
\draw [thick]  (9,3.5) arc (-90:90:10mm);
\draw [thick] (v1) -- (v2) -- (v4);
\draw [thick] (v5) -- (v4);
\node at (7.75,3.75) {$v_{4}'$};
\node [blue] at (8.5,3.5) {$e_1$};
\node [blue] at (8.75,4.5) {$e_3$};

\node at (11,4.5) {$+$};
\node [place, label=above:$v_{b+4}$] at (-1,0.5) {};
\node at (-0.5,0.5) {$\cdots$};
\node [place, label=above:$v_5$] (v12) at (0,0.5) {};
\node [sample, label=below:$v_{4}$] (v13) at (1,1) {};
\node [place, label=below:$v_{3}$] (v16) at (2,1.5) {};
\node [place, label=left:$v_{2}$] at (3,0.5) {};
\node [place, label=above:$v_{1}$] (v15) at (2,-0.5) {};
\draw [thick]  (2,-0.5) arc (-90:90:10mm);
\node [sample, label=above:$v_{4}'$] (v14) at (1,0) {};
\draw [thick] (v14) -- (v13);
\draw [thick] (v12) -- (v13) -- (v15);
\node [blue] at (1.75,0.5) {$e_3$};
\node [blue] at (0.75,0.5) {$e_2$};

\node at (4.5,0.5) {$\sim$};
\node [place, label=above:$v_{b+4}$] at (6,0.5) {};
\node at (6.5,0.5) {$\cdots$};
\node [place, label=above:$v_5$] (v17) at (7,0.5) {};
\node [sample, label=below:$v_{4}$] (v18) at (8,1) {};
\node [place, label=below:$v_{3}$] (v19) at (9,1.5) {};
\node [place, label=left:$v_{2}$] at (10,0.5) {};
\node [place, label=above:$v_{1}$] (v20) at (9,-0.5) {};
\draw [thick] (9,-0.5) node (v22) {} arc (-90:90:10mm);
\node [sample, label=above:$v_{4}'$] (v21) at (8,0) {};
\draw [thick] (v17) -- (v18);
\draw [thick] (v20) edge (v18);
\node [blue] at (8.75,0.5) {$e_3$};
\draw [thick] (v18) -- (v21) -- (v20);
\node [blue] at (7.75,0.5) {$e_2$};
\node [blue] at (8.5,-0.5) {$e_1$};
\end{tikzpicture}
\caption{
%An illustration of the Lemma \ref{PC'-2}
Decomposition of $H_b-\{v_3v_4,\,v_3v_4'\}$.}\label{gra-h5}
\end{figure}
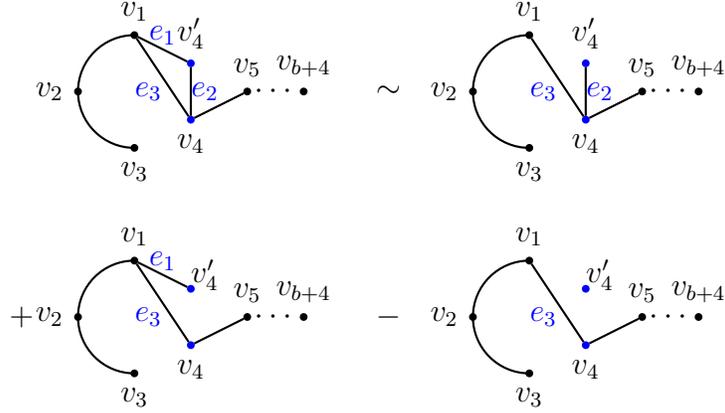

Using the same method, we deduce that
\begin{align}\label{equ-H-4}
X_{H_b-\{v_3v_4,v_3v_4'\}}=X_{P_{b+4,4}}+X_{P_{b+4,3}}-e_1X_{P_{b+4}},
\end{align}
as shown in Figure \ref{gra-h5}.

Substituting \eqref{equ-H-2}, \eqref{equ-H-3} and \eqref{equ-H-4} into \eqref{equ-H-1}, we get
\begin{align}\label{xhb-pmn}
X_{H_b} =& 2X_{T^{\langle 4,b+1 \rangle}}+2X_{T^{\langle 5,b \rangle}}+e_1(X_{T^{\langle 4,b \rangle}}+X_{P_{b+4}})-X_{P_{b+1}}X_{C_4}-2X_{P_{b+4,4}} \notag \\
&-X_{P_{b+4,3}}-X_{P_{b+5}}.
\end{align}

%For $T^{\langle a,b \rangle}$ and $P_{m,n}$ we have the following expansions
%\begin{align}\label{pmn}
%X_{T^{\langle a,b \rangle}} =&
%%\sum_{i=0}^b X_{P_{i}}X_{C_{a+b-i}} - b X_{P_{a+b}} \quad (\text{here we set } X_{P_0} = 1) \\
%(a-1)X_{P_{a+b}} - \sum_{i=2}^{a-1}X_{C_i}X_{P_{a+b-i}}. \\ \label{tab}
%X_{P_{m,n}} =& X_{P_{m+1}}+e_1 X_{P_m}-X_{P_n}X_{P_{m-n+1}} .
%\end{align}
%We can prove the above results by applying Proposition \ref{dec-thm} and Corollary \ref{dec-cor} in the same manner as before, and the details are omitted here.  ********
Now it remains to calculate $X_{T^{\langle a,b \rangle}}$ and $X_{P_{m,n}}$. For $P_{m,n}$, we label $v_{n-1}v_n$, $v_nv_{m+1}$ by $e_1$, $e_2$ respectively, as shown in Figure \ref{gra-pmn}. By using Proposition \ref{dec-cor}, we obtain
\begin{align}\label{pmn}
X_{P_{m,n}} =X_{P_{m+1}}+e_1 X_{P_m}-X_{P_n}X_{P_{m-n+1}}.
\end{align}
For $T^{\langle a,b \rangle}$, we label $v_1v_a$, $v_{a-1}v_a$ by $e_1$, $e_2$ respectively, as shown in Figure \ref{gra-tab}. By using Proposition \ref{dec-cor}, we have
\begin{align}\label{tab1}
X_{T^{\langle a,b \rangle}}=X_{T^{\langle a-1,b+1 \rangle}}+X_{P_{a+b}}-X_{C_{a-1}}X_{P_{b+1}}.
\end{align}
Iteration of \eqref{tab1} leads to
\[X_{T^{\langle a,b \rangle}}=X_{T^{\langle 3,a+b-3 \rangle}}+(a-3)X_{P_{a+b}}-\sum_{i=3}^{a-1}X_{C_i}X_{P_{a+b-i}}.\]
On the other hand, by Proposition \ref{dec-thm}, we have
\[X_{T^{\langle 3,a+b-3 \rangle}}=2X_{P_{a+b}}-X_{C_2}X_{P_{a+b-2}},\]
as shown in Figure \ref{gra-tab1}. Hence we deduce that
\begin{align}\label{tab}
X_{T^{\langle a,b \rangle}} =&
%\sum_{i=0}^b X_{P_{i}}X_{C_{a+b-i}} - b X_{P_{a+b}} \quad (\text{here we set } X_{P_0} = 1) \\
(a-1)X_{P_{a+b}} - \sum_{i=2}^{a-1}X_{C_i}X_{P_{a+b-i}}.
\end{align}

\begin{figure}[ht]
\centering
\begin{tikzpicture}
[place/.style={thick,fill=black!100,circle,inner sep=0pt,minimum size=1mm,draw=black!100},rotate=180,scale=0.75]
\node [place, label=above:$v_{a+b}$] (v3) at (-1,7) {};
\node [place, label=above:$v_{a+1}$] (v4) at (0,7) {};
\node [place, label=left:$v_{a}$] (v5) at (1,7) {};
\node at (-0.5,7) {$\cdots$};
\draw [thick] (2,6) arc (-90:-270:10mm);
\draw [thick] (2,6) [densely dashed] arc (-90:90:10mm);
\node [place, label=above:$v_{1}$] at (2,6) {};
\node [place, label=below:$v_{a-1}$] at (2,8) {};
\node [blue] at (1.25,6) {$e_1$};
\node [blue] at (1.25,8) {$e_2$};
\draw [thick] (v5) edge (v4);

\node at (-2,7) {$\sim$};
\node [place, label=above:$v_{1}$] (v1) at (-4,6) {};
\node [place, label=left:$v_{a}$] (v7) at (-5,7) {};
\node [place, label=below:$v_{a-1}$] (v2) at (-4,8) {};
\draw [thick] (-4,8) arc (90:180:10mm);
\draw [thick] (-4,8) [densely dashed] arc (90:-90:10mm);
\node [blue] at (-4.75,8) {$e_2$};
\draw [thick] (v1) -- (v2);
\node [blue] at (-3.75,7) {$e_3$};
\node [place, label=above:$v_{a+1}$] (v6) at (-6,7) {};
\node at (-6.5,7) {$\cdots$};
\node [place, label=above:$v_{a+b}$] at (-7,7) {};
\draw [thick] (v6) -- (v7);

\node [place, label=above:$v_{1}$] at (2,9.5) {};
\node [place, label=below:$v_{a-1}$] at (2,11.5) {};
\node [place, label=left:$v_{a}$] (v8) at (1,10.5) {};
\node [place, label=above:$v_{a+1}$] (v9) at (0,10.5) {};
\node at (-0.5,10.5) {$\cdots$};
\node [place, label=above:$v_{a+b}$] at (-1,10.5) {};
\draw [thick] (2,9.5) arc (-90:-180:10mm);
\draw [thick] (2,9.5) [densely dashed] arc (-90:90:10mm);
\draw [thick] (v8) -- (v9);
\node [blue] at (1.25,9.5) {$e_1$};

\node [place, label=above:$v_{1}$] (v10) at (-4,9.5) {};
\node [place, label=left:$v_{a}$] (v12) at (-5,10.5) {};
\node [place, label=below:$v_{a-1}$] (v11) at (-4,11.5) {};
\draw [thick] (-4,11.5) [densely dashed] arc (90:-90:10mm);
\draw [thick] (v10) -- (v11);
\node [blue] at (-3.75,10.5) {$e_3$};
\node [place, label=above:$v_{a+1}$] (v13) at (-6,10.5) {};
\node at (-6.5,10.5) {$\cdots$};
\node [place, label=above:$v_{a+b}$] at (-7,10.5) {};
\draw [thick] (v13) -- (v12);
\node at (3.5,10.5) {$+$};
\node at (-2,10.5) {$-$};
\end{tikzpicture}
\caption{Decomposition of $T^{\langle a,b \rangle}$.}\label{gra-tab}
\end{figure}

\newpage

\begin{figure}[h]
\centering
\begin{tikzpicture}
[place/.style={thick,fill=black!100,circle,inner sep=0pt,minimum size=1mm,draw=black!100},rotate=180]
\node [place, label=above:$v_{a+b}$] (v3) at (-1,7) {};
\node [place, label=above:$v_{4}$] (v4) at (0,7) {};
\node [place, label=above:$v_{3}$] (v5) at (1,7) {};
\node at (-0.5,7) {$\cdots$};
\node [place, label=above:$v_{1}$] (v14) at (2,6.5) {};
\node [place, label=below:$v_{2}$] (v15) at (2,7.5) {};
\node [blue] at (1.5,6.5) {$e_1$};
\node [blue] at (1.5,7.5) {$e_2$};
\draw [thick] (v5) edge (v4);
\draw [thick] (v5) -- (v14) -- (v15) -- (v5);
\node [blue] at (2.25,7) {$e_3$};

\node at (-2,7) {$\sim$};
\node [place, label=above:$v_{1}$] (v1) at (-3,6.5) {};
\node [place, label=above:$v_{3}$] (v7) at (-4,7) {};
\node [place, label=below:$v_{2}$] (v2) at (-3,7.5) {};
\node [blue] at (-3.5,7.5) {$e_2$};
\node [blue] at (-2.75,7) {$e_3$};
\node [place, label=above:$v_{4}$] (v6) at (-5,7) {};
\node at (-5.5,7) {$\cdots$};
\node [place, label=above:$v_{a+b}$] at (-6,7) {};
\draw [thick] (v6) -- (v7);
\draw [thick] (v1) -- (v2) -- (v7);

\node at (2.75,10.5) {$+$};
\node [place, label=above:$v_{1}$] at (2,10) {};
\node [place, label=below:$v_{2}$] (v16) at (2,11) {};
\node [place, label=above:$v_{3}$] (v8) at (1,10.5) {};
\node [place, label=above:$v_{4}$] (v9) at (0,10.5) {};
\node at (-0.5,10.5) {$\cdots$};
\node [place, label=above:$v_{a+b}$] at (-1,10.5) {};
\draw [thick] (v8) -- (v9);
\node [blue] at (1.5,10) {$e_1$};
\node [blue] at (2.25,10.5) {$e_3$};
\draw [thick] (v8) -- (2,10) -- (v16);

\node at (-2,10.5) {$-$};
\node [place, label=above:$v_{1}$] (v10) at (-3,10) {};
\node [place, label=above:$v_{3}$] (v12) at (-4,10.5) {};
\node [place, label=below:$v_{2}$] (v11) at (-3,11) {};
\draw [thick] (v10) -- (v11);
\node [blue] at (-2.75,10.5) {$e_3$};
\node [place, label=above:$v_{a+1}$] (v13) at (-5,10.5) {};
\node at (-5.5,10.5) {$\cdots$};
\node [place, label=above:$v_{a+b}$] at (-6,10.5) {};
\draw [thick] (v13) -- (v12);
\end{tikzpicture}
\caption{Decomposition of $T^{\langle 3,\,a+b-3 \rangle}$.}\label{gra-tab1}
\end{figure}

Finally, putting together \eqref{xhb-pmn}, \eqref{pmn} and \eqref{tab}, we get
\begin{align}\label{xhb}
X_{H_b} =& 10 X_{P_{b+5}} + e_1X_{P_{b+4}} - 4X_{C_2}X_{P_{b+3}} - (4X_{C_3}+e_1X_{C_2}-X_{P_3})X_{P_{b+2}} \notag \\
&- (3X_{C_4} + e_1X_{C_3} - 2X_{P_4}) X_{P_{b+1}}.
\end{align}
Combining \eqref{xt1}, \eqref{xs}, \eqref{xhb} and the explicit formulas
\begin{align*}
X_{C_4} &= 2e_{(2,2)} + 12e_{4}, \quad X_{C_3} = 6e_{3}, \quad X_{C_2} = 2e_{2} \\
X_{P_4} &= 2e_{(2,2)} + 2e_{(3,1)} + 4e_4, \quad X_{P_3} = e_{(2,1)} + 3e_3.
\end{align*}
we obtain the desired expansion \eqref{t'}. \qed

To prove the non-$e$-positivity of $X^{\langle 4,b \rangle}_{v_4}$, we need to calculate the coefficients in the expansion of $X_{P_n}$ in terms of elementary symmetric functions. Based on the generating function of $X_{P_n}$ given by Stanley \cite{Sta95}, Wolfe showed the following explicit formula.

\begin{lem}\cite[Theorem 3.2]{Wol98}\label{pathcoe}
Let $\lambda = \langle 1^{r_1},2^{r_2}, \ldots, n^{r_n} \rangle$ be a partition of $n$. Then the coefficient $c_{\lambda}$ of $e_{\lambda}$ in the expansion of $X_{P_n}$ is given by
\[
	c_{\lambda} = \binom{r_1+\cdots+r_n}{r_1,\ldots,r_n}\prod_{j=1}^{n}(j-1)^{r_j} + \sum_{i=1}^n \left(\binom{(r_1+\cdots+r_n) - 1}{r_1,\ldots,r_i-1,\ldots,r_n} (i-1)^{r_i-1} \prod_{\substack{j = 1 \\ j \neq i}}^{n} (j-1)^{r_j} \right),
\]
where by convention $d^0 = 1$ for all $d \ge 0$ and $\binom{r_1+\cdots+r_n}{r_1,\ldots,r_n} = 0$ if some $r_i < 0$.
\end{lem}

We proceed to show that certain coefficient is negative in the expansion of $X_{T^{\langle 4,b \rangle}_{v_4}}$ in terms of elementary symmetric functions, which together with Theorem \ref{tadepos} gives a series of counterexamples to Conjecture \ref{mconj}.

\begin{thm}
For any $b \geq 1$, $X_{T^{\langle 4,b \rangle}_{v_4}}$ is not $e$-positive. In particular, for $k \ge 0$, we have
\begin{align*}
[e_{(3^{k+2})}]X_{T^{\langle 4, 3k+1 \rangle}_{v_4}} =& -3\cdot 2^{k+1}, \\[5pt]
[e_{(3^{k+2},1)}]X_{T^{\langle 4, 3k+2 \rangle}_{v_4}} =& - 2^{k+2}, \\[5pt]
[e_{(3^{k+2},2)}]X_{T^{\langle 4, 3k+3 \rangle}_{v_4}} =& -(3k+1)\cdot 2^{k+1}.
\end{align*}
\end{thm}

\pf By the definition of elementary symmetric functions, if $b = 3k+1$, then only the terms $20 X_{P_{b+5}}$ and $-42 e_3 X_{P_{b+2}}$ in the expansion \eqref{t'} contribute to the coefficient $[e_{(3^{k+2})}]X_{T^{\langle 4, b \rangle}_{v_4}}$. Precisely, we have
\begin{align*}
	[e_{(3^{k+2})}]X_{T^{\langle 4, 3k+1 \rangle}_{v_4}} = 20 [e_{(3^{k+2})}]X_{P_{3k+6}} - 42 [e_{(3^{k+1})}] X_{P_{3k+3}}.
\end{align*}
By Lemma \ref{pathcoe},
\begin{align*}
	[e_{(3^{k+2})}]X_{P_{3k+6}} =& \binom{k+2}{k+2}2^{k+2} + \binom{k+1}{k+1}2^{k+1} = 6 \cdot 2^{k}, \\[5pt]
	[e_{(3^{k+1})}] X_{P_{3k+3}} =& \binom{k+1}{k+1}2^{k+1} + \binom{k}{k}2^{k} = 3 \cdot 2^{k},
\end{align*}
and hence
\[
	[e_{(3^{k+2})}]X_{T^{\langle 4, 3k+1 \rangle}_{v_4}} = 20 \cdot 6 \cdot 2^{k} - 42 \cdot 3 \cdot 2^{k} = -3\cdot 2^{k+1}.
\]

For the other two cases, it is straightforward to verify that
\begin{align*}
	[e_{(3^{k+2},1)}]X_{T^{\langle 4, 3k+2 \rangle}_{v_4}} =& 20[e_{(3^{k+2},1)}]X_{P_{3k+7}} + 2[e_{(3^{k+2})}]X_{P_{3k+6}} - 42[e_{(3^{k+1},1)}]X_{P_{3k+4}} \\
	&- 4[e_{(3^{k+1})}]X_{P_{3k+3}} \\[5pt]
	=& 20 \cdot 2^{k+2} + 2 \cdot (6 \cdot 2^{k}) - 42 \cdot 2^{k+1} - 4 \cdot (3 \cdot 2^{k}) \\[5pt]
	=& -2^{k+2}
\end{align*}
and
\begin{align*}
	[e_{(3^{k+2},2)}]X_{T^{\langle 4, 3k+3 \rangle}_{v_4}} =& 20[e_{(3^{k+2},2)}]X_{P_{3k+8}} - 16[e_{(3^{k+2})}]X_{P_{3k+6}} - 42[e_{(3^{k+1},2)}]X_{P_{3k+5}} \\
	&- 4[e_{(3^{k+1})}]X_{P_{3k+3}} \\[5pt]
	=& 20 ( (k+3) 2^{k+2} + 2^{k+2} + (k+2) 2^{k+1} ) - 16 \cdot (6 \cdot 2^{k}) \\
	&- 42 ( (k+2) 2^{k+1} + 2^{k+1} + (k+1) 2^k ) - 4 \cdot (3 \cdot 2^{k}) \\[5pt]
	=& -(3k+1)\cdot 2^{k+1}.
\end{align*}
This completes the proof.
\qed

\section{Counterexample to Conjollary \ref{mconj-coro}} \label{cx-s-pos}

\begin{figure}[ht]
\centering
\begin{tikzpicture}[scale = 1]
    \fill (0,0) circle (0.3ex);
    \fill (1,0) circle (0.3ex);
    \fill[red] (2,0) circle (0.3ex);
    \node[right] (w) at (2.2,0) {$w$};
    \node[right] (x) at (3.2,1) {$x$};
    \node[right] (y) at (3.2,-1) {$y$};
    \node[below] (u) at (0,-0.2) {$u$};
    \node[below] (v) at (1,-0.2) {$v$};
    \node[below] (f) at (1.5,-1.5) {$F$};
    \node[below] (f1) at (6.5,-1.5) {$F^{(1)}_w$};
    \node[below] (f) at (11.5,-1.5) {$F^{(2)}_w$};
    \fill (3,1) circle (0.3ex);
    \fill (3,-1) circle (0.3ex);
    \draw (0,0) -- (1,0);
    \draw (1,0) -- (2,0);
    \draw (2,0) -- (3,1);
    \draw (2,0) -- (3,-1);
    \fill (5,0) circle (0.3ex);
    \fill (6,0) circle (0.3ex);
    \fill[red] (7,0.5) circle (0.3ex);
    \fill[red] (7,-0.5) circle (0.3ex);
    \fill (8,1) circle (0.3ex);
    \fill (8,-1) circle (0.3ex);
    \draw (5,0) -- (6,0);
    \draw (6,0) -- (7,0.5);
    \draw (6,0) -- (7,-0.5);
    \draw (7,0.5) -- (8,1);
    \draw (7,-0.5) -- (8,-1);
    \draw (7,0.5) -- (8,-1);
    \draw (7,-0.5) -- (8,1);
    \draw (7,0.5) -- (7,-0.5);
    \fill (10,0) circle (0.3ex);
    \fill (11,0) circle (0.3ex);
    \fill[red] (12,0.5) circle (0.3ex);
    \fill[red] (12,-0.5) circle (0.3ex);
    \fill[red] (13,0) circle (0.3ex);
    \fill (13.5,1) circle (0.3ex);
    \fill (13.5,-1) circle (0.3ex);
    \draw (10,0) -- (11,0);
    \draw (11,0) -- (12,0.5);
    \draw (11,0) -- (12,-0.5);
    \draw (11,0) -- (13,0);
    \draw (12,0.5) -- (12,-0.5);
    \draw (12,-0.5) -- (13,0);
    \draw (13,0) -- (12,0.5);
    \draw (12,0.5) -- (13.5,1);
    \draw (12,0.5) -- (13.5,-1);
    \draw (12,-0.5) -- (13.5,1);
    \draw (12,-0.5) -- (13.5,-1);
    \draw (13,0) -- (13.5,1);
    \draw (13,0) -- (13.5,-1);
\end{tikzpicture}
    \caption{The fork graph and two of its expansions.}\label{fork}
\end{figure}
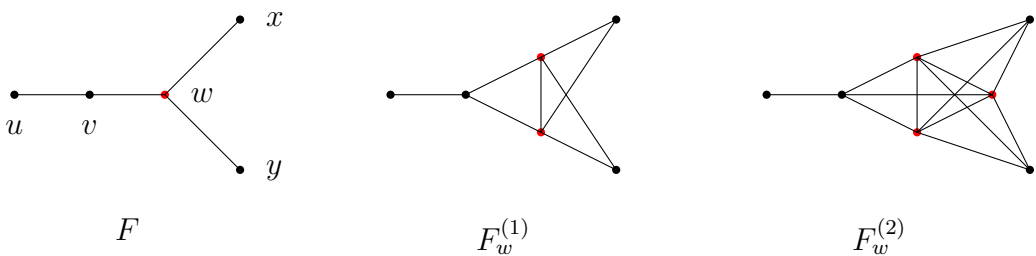

Although Conjecture \ref{mconj} was disproved in the last section by showing that all $X_{T^{\langle 4, b \rangle}_{v_4}}$ are not $e$-positive, experimental results show that they are still $s$-positive. This means that Conjollary \ref{mconj-coro} is possibly true. The main objective of this section is to provide a counterexample to Conjollary \ref{mconj-coro}.
For an $e$-positive graph $G$ and a vertex $v\in V(G)$, instead of just considering $G_v$, we shall study the $s$-positivity of general $G^{(k)}_v$ for any fixed $k \ge 1$, where $G^{(k)}_v$ is  the clan graph obtained from $G$ by replacing $v$ with a complete graph $K_{k+1}$, and in particular $G_v = G^{(1)}_v$. Here we take $G$ to be the fork graph $F$, whose vertices are labeled as in Figure \ref{fork}, and focus on the $s$-positivity of $F^{(k)}_w$.

With the help of Sage \cite{Sage}, we find that
\begin{align}
    X_F = 5e_{(5)} + 7e_{(4,1)} + e_{(3,2)} + 2e_{(3,1,1)} + e_{(2,2,1)}, \label{XF-1}
\end{align}
which is $e$-positive.
The following result, together with \eqref{XF-1}, implies that Conjollary \ref{mconj-coro} fails for the fork graph $F$.

\begin{thm}\label{main-thm-sect3}
Let $F$ and $w$ be given as in Figure \ref{fork}. Then for any $k\geq 1$
\begin{align}\label{XF}
	X_{F^{(k)}_w}=& (k + 1)! s_{(3,2,1^k)} + {(k + 3)} (k + 1)! s_{(3,1^{k+2})} -2(k + 1)! s_{(2^3,1^{k-1})}  \notag \\
	&+ 5{(k + 1)}(k + 1)! s_{(2^2,1^{k+1})} + {(k^{3} + 5 k^{2} + 14 k + 14)} (k + 1)! s_{(2,1^{k+3})} \notag
\\
  &+ 2 {(k + 2)}^{3} (k + 1)! s_{(1^{k+5})}.
\end{align}
In particular, $X_{F^{(k)}_w}$ is not $s$-positive.
\end{thm}

\pf We first use Proposition \ref{stanley-lemma} to give the monomial expansion of $X_{F^{(k)}_w}$.
A little thought shows that any stable partition of $F^{(k)}_w$ must be of type
$(3,2,1^k)$, $(3,1^{k+2})$, $(2^2,1^{k+1})$, $(2,1^{k+3})$ or $(1^{k+5})$.
It remains to determine the number $a_{\lambda}$ of stable partitions of each aforementioned type $\lambda$.
We take $\lambda=(2^2,1^{k+1})$ as an example to give an illustration. To compute $a_{(2^2,1^{k+1})}$ we keep the labels $u,v,x,y$ of those four vertices of $F$ in $F^{(k)}_w$ and label the remaining vertices of $F^{(k)}_w$ by $w_{0},w_{1},\ldots,w_{k}$.
Then every stable partition of type $(2^2,1^{k+1})$ is of the form $\{u,w_i\}/\{v,x\}/\cdots$, $\{u,w_i\}/\{v,y\}/\cdots$, $\{u,w_i\}/\{x,y\}/\cdots$, $\{u,x\}/\{v,y\}/\cdots$, or $\{u,y\}/\{v,x\}/\cdots$, where $0 \le i \le k$. Hence $a_{2^2,1^{k+1}} = 3(k+1)+2 = 3k+5$. The number of stable partitions of other types may be obtained in the same way. Finally, we get
\begin{align}\label{eq-am}
    X_{F^{(k)}_w} = (k+1)\tilde{m}_{(3,2,1^k)} + 2\tilde{m}_{(3,1^{k+2})} + (3k+5)\tilde{m}_{(2^2,1^{k+1})} + (k+6)\tilde{m}_{(2,1^{k+3})} + \tilde{m}_{(1^{k+5})}.
\end{align}

We proceed to give the Schur function expansion of $X_{F^{(k)}_w}$. At first, one should note that by Proposition \ref{Kosnum}, if we order the bases in a linear order compatible with dominance order, then the transition matrices between the monomial symmetric functions $m_{\lambda}$ and the Schur functions $s_{\lambda}$ are triangular and hence invertible. Precisely, to express the monomial symmetric functions in $X_{F^{(k)}_w}$ by Schur functions, we only need to consider the $s_{\lambda}$'s with $\lambda \le \mu$ for some $\mu$ in $\{(3,2,1^k),(3,1^{k+2}),(2^2,1^{k+1}),(2,1^{k+3}),(1^{k+5})\}$ under dominance order. For any $k \ge 1$, we list all these Schur functions and the transition matrix as follows, and the computation of the entries in the matrix will be discussed in the next paragraph.
{\footnotesize
\begin{align*}
    \begin{pmatrix}
    s_{(3,2,1^k)} \\
    s_{(3,1^{k+2})} \\
    s_{(2^3,1^{k-1})} \\
    s_{(2^2,1^{k+1})} \\
    s_{(2,1^{k+3})} \\
    s_{(1^{k+5})}
    \end{pmatrix}
=\begin{pmatrix}
    1 & k + 1 & 2 & 2  k + 2 & {(k + 3)} {(k + 1)} & {(k + 5)} {(k + 3)} {(k + 1)} / 3\\
    0 & 1 & 0 & 1 & k + 3 & {(k + 4)} {(k + 3)}/2 \\
    0 & 0 & 1 & k & {(k + 3)} k/2 & {(k + 5)} {(k + 4)} k/6 \\
    0 & 0 & 0 & 1 & k + 2 & {(k + 5)} {(k + 2)}/2 \\
    0 & 0 & 0 & 0 & 1 & k + 4 \\
    0 & 0 & 0 & 0 & 0 & 1
    \end{pmatrix}
        \begin{pmatrix}
    m_{(3,2,1^k)} \\
    m_{(3,1^{k+2})} \\
    m_{(2^3,1^{k-1})} \\
    m_{(2^2,1^{k+1})} \\
    m_{(2,1^{k+3})} \\
    m_{(1^{k+5})}
    \end{pmatrix}.
\end{align*}
}

\begin{figure}[ht]
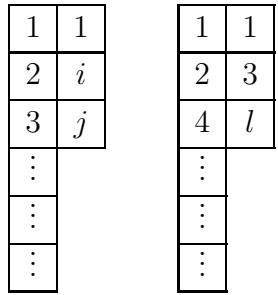

\begin{center}
\begin{ytableau}[]
1 & 1 \\
2 & i \\
3 & j \\
\vdots & \none \\
\vdots & \none \\
\vdots & \none
\end{ytableau} \qquad
\begin{ytableau}[]
1 & 1 \\
2 & 3 \\
4 & l \\
\vdots & \none \\
\vdots & \none \\
\vdots & \none
\end{ytableau}
\end{center}
\caption{The possible Young tableaux of $K_{(2^3,1^{k-1}),(2,1^{k+3})}$.}\label{Kos-ex}
\end{figure}

As stated in Proposition \ref{smes}, the coefficients in the above transition matrix are Kostka numbers. The last column, whose entries are all of the form $K_{\lambda,(1^{k+5})}$, are easily computed by Proposition \ref{hlf}, the hook length formula. For the other numbers, it is also not difficult to give explicit formulas by using the definition of Kostka numbers. To illustrate, we take $K_{(2^3,1^{k-1}),(2,1^{k+3})}$ as an example. To fill a Young diagram of shape $(2^3,1^{k-1})$ with two 1's, one 2, one 3, \ldots, and one $k+4$, we must put all 1's into the first row and the number 2 into the first box of the second row. Then we may put 3 into the second box of the second row or the first box of the third row, and in the latter case the number 4 must be filled into the first box of the third row. We draw these two cases in Figure \ref{Kos-ex}.
In the left diagram, if we fix $i$ and $j$, then there is only one way to fill the other numbers, and the same thing happens if we fix $l$ in the right diagram. Then since $i < j$ and $l$ may be chosen arbitrarily, the total number of the possible Young tableaux is $\binom{k+1}{2} + k = (k+3)k/2$. We omit the computation of the other Kostka numbers in the above matrix.

From \eqref{eq-am} it follows that 
{\footnotesize
\begin{align*}
    X_{F^{(k)}_w}
    =& ((k+1)!, 2(k+2)!, 0, 2(3k+5)(k+1)!, (k+6)(k+3)!, (k+5)!) \cdot		
    \begin{pmatrix}
    m_{(3,2,1^k)} \\
    m_{(3,1^{k+2})} \\
    m_{(2^3,1^{k-1})} \\
    m_{(2^2,1^{k+1})} \\
    m_{(2,1^{k+3})} \\
    m_{(1^{k+5})}
    \end{pmatrix}.
\end{align*}}

Putting together the above two identities, it is routine to verify that
{\footnotesize
\begin{align*}
    X_{F^{(k)}_w}
    =& ((k+1)!, 2(k+2)!, 0, 2(3k+5)(k+1)!, (k+6)(k+3)!, (k+5)!)  \cdot
    \notag \\
    & \quad \begin{pmatrix}
    1 & k + 1 & 2 & 2  k + 2 & {(k + 3)} {(k + 1)} & {(k + 5)} {(k + 3)} {(k + 1)} / 3\\
    0 & 1 & 0 & 1 & k + 3 & {(k + 4)} {(k + 3)}/2 \\
    0 & 0 & 1 & k & {(k + 3)} k/2 & {(k + 5)} {(k + 4)} k/6 \\
    0 & 0 & 0 & 1 & k + 2 & {(k + 5)} {(k + 2)}/2 \\
    0 & 0 & 0 & 0 & 1 & k + 4 \\
    0 & 0 & 0 & 0 & 0 & 1
    \end{pmatrix}^{-1}
    \begin{pmatrix}
    s_{(3,2,1^k)} \\
    s_{(3,1^{k+2})} \\
    s_{(2^3,1^{k-1})} \\
    s_{(2^2,1^{k+1})} \\
    s_{(2,1^{k+3})} \\
    s_{(1^{k+5})}
    \end{pmatrix}
    \notag \\[8pt]
    =& ((k+1)!, 2(k+2)!, 0, 2(3k+5)(k+1)!, (k+6)(k+3)!, (k+5)!)  \cdot \notag \\
    & \quad
    \begin{pmatrix}
    1 & -k - 1 & -2 & k - 1 & 2k+2 & -(k+2)(k+1) \\
    0 & 1 & 0 & -1 & -1 & k+3 \\
    0 & 0 & 1 & -k & (k+1)k/2 & -(k + 2)(k + 1)k/6 \\
    0 & 0 & 0 & 1 & -k - 2  & (k+3)(k+2)/2 \\
    0 & 0 & 0 & 0 & 1 & -k - 4 \\
    0 & 0 & 0 & 0 & 0 & 1
    \end{pmatrix}
    \begin{pmatrix}
    s_{(3,2,1^k)} \\
    s_{(3,1^{k+2})} \\
    s_{(2^3,1^{k-1})} \\
    s_{(2^2,1^{k+1})} \\
    s_{(2,1^{k+3})} \\
    s_{(1^{k+5})}
    \end{pmatrix}
    \notag \\[8pt]
    =& (k + 1)! s_{(3,2,1^k)} + {(k + 3)} (k + 1)! s_{(3,1^{k+2})} -2(k + 1)! s_{(2^3,1^{k-1})} + 5{(k + 1)}(k + 1)! s_{(2^2,1^{k+1})} \notag \\
	& + {(k^{3} + 5 k^{2} + 14 k + 14)} (k + 1)! s_{(2,1^{k+3})} + 2 {(k + 2)}^{3} (k + 1)! s_{(1^{k+5})},
\end{align*}}
as desired. This completes the proof.  \qed

\begin{figure}[htbp]
\begin{center}
\begin{tikzpicture}[scale=0.5]
\tyng(0cm, 0cm, 2^3,1^1)
\tyng(0cm, 6cm, 1)
\coordinate (L1) at (0.5, 0.5);
\coordinate (L2) at (0.5, 1.5); 
\coordinate (L3) at (0.5, 2.5); 
\coordinate (L4) at (0.5, 3.5);
\coordinate (L5) at (0.5, 6.5);  
\coordinate (R2) at (1.5, 1.5);
\coordinate (R1) at (1.5, 0.5); 
\coordinate (a)  at (0.5, 5.3);
\coordinate (b)  at (0.5, 5); 
\coordinate (c)  at (0.5, 4.7); 
\draw[RimHook]
(L3) -- (1.5, 2.5) -- (R2);
\draw[RimHook]
(L1) -- (R1);
\draw[OuterBoundary]
(0,0)
-- (2, 0)
-- (2, 3)
-- (1, 3)
-- (1, 7)
-- (0, 7)
-- (0, 0);
\foreach \s in {L1,L2, L3, L4, L5, R2, R1}
  \shade[RHiball](\s) circle(.1);
\foreach \t in {a, b, c}
  \shade[apball](\t) circle(.03);
\begin{scope}[xshift=6cm]
\tyng(0cm, 0cm, 2^3,1^1)
\tyng(0cm, 6cm, 1)
\coordinate (L1) at (0.5, 0.5);
\coordinate (L2) at (0.5, 1.5); 
\coordinate (L3) at (0.5, 2.5); 
\coordinate (L4) at (0.5, 3.5);
\coordinate (L5) at (0.5, 6.5);  
\coordinate (R3) at (1.5, 2.5);
\coordinate (R2) at (1.5, 1.5);
\coordinate (R1) at (1.5, 0.5); 
\coordinate (a)  at (0.5, 5.3);
\coordinate (b)  at (0.5, 5); 
\coordinate (c)  at (0.5, 4.7); 
\draw[RimHook]
(L3) -- (R3);
\draw[RimHook]
(L2) -- (R2) -- (R1);
\draw[OuterBoundary]
(0,0)
-- (2, 0)
-- (2, 3)
-- (1, 3)
-- (1, 7)
-- (0, 7)
-- (0, 0);
\foreach \s in {L1,L2, L3, L4, L5, R3, R1}
  \shade[RHiball](\s) circle(.1);
\foreach \t in {a, b, c}
  \shade[apball](\t) circle(.03);
\end{scope}
\end{tikzpicture}
\end{center}
\caption{The two special rim hook tabloids in $\mathcal{T}_{(2^3,\,1^{k-1})}$ for the graph $F^{(k)}_w$.}\label{fig:tabloids}
\end{figure}
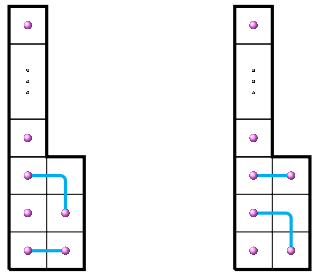

We remark that all the coefficients in \eqref{XF} are combinatorially computable with the aid of special rim hook tabloids, see the third author and Wang~\cite[Theorem~3.1]{WW20}. Confirming the negativity of the coefficient $\brk[s]1{s_{(2^3,\,1^{k-1})}}X_{F^{(k)}_w}$ is speticularly rather easy.
Explicitly speaking, the coefficient is
\begin{equation}\label{eq:tabloid}
\brk[s]1{s_{\lambda}}X_{G}
=\sum_{T\in\mathcal{T}_{\lambda}}(-1)^{\abs{W_T}}N_T,
\end{equation}
where $\lambda=(2^3,\,1^{k-1})$, $G=F^{(k)}_w$,
$\mathcal{T}_{\lambda}$ is the set of special rim hook tabloids~$T$ of shape~$\lambda$ such that $G$ contains a stable partition of type $\kappa_T$ (the partition consisting of the rim hook lengths in $T$),
$W_T$ is the set of rim hooks of $T$ that span an even number of rows, and $N_T$ is the number of stable partitions of $G$ of type $\kappa_T$ whose blocks of the same size are ordered. Because of the special structure of $G$, 
there are only six types $\kappa_T$ of which~$G$ contains a stable partition:
\[
\kappa=(3,2,1^k),\quad
(3,1^{k+2}),\quad
(2^3,1^{k-1}),\quad
(2^2,1^{k+1}),\quad
(2,1^{k+3}),\quad\text{and}\ 
(1^{k+5}).
\]
Observe that there does not exist a special rim hook tabloid of shape $\lambda$ and any of the above types except $\kappa$; there are only two such tabloids of type $\kappa$ as shown in Figure~\ref{fig:tabloids}.
From definition, it is clear that $\abs{W_T}=1$ for each of the two tabloids $T\in\mathcal{T}_\lambda$, 
equally clear is the existence of a stable partition of type $\kappa$.
Therefore, every summand in \eqref{eq:tabloid} is negative and the desired negativity follows. 
%Interested reader can verify that \eqref{eq:tabloid} has only two summands of the same value $-(k+1)!$.

\section{Further directions}\label{sect-last}

Although Conjecture \ref{mconj} and Conjollary \ref{mconj-coro} are not valid in general, its restrictions to certain special cases are expected to be true.  

The first special case relates to unit interval graphs. In 2013, Guay-Paquet \cite{GP13} reduced Stanley's $\mathbf{(3+1)}$-free conjecture to proving the $e$-positivity of unit interval graphs. Recall that a graph $G$ is said to be a \textit{unit interval graph} if we can find a multiset $S$ of unit intervals (or just intervals of the same length) of the real line  such that there is a bijection $f: V(G) \to S$ satisfying the condition that two vertices $u$ and $v$ are adjacent in $G$ if and only if the two intervals $f(u)$ and $f(v)$ intersect. If $S$ is a multiset of unit intervals and $I \in S$, then denote by $S^{(k)}_I$ the multiset of unit intervals obtained from $S$ by adding $k$ copies of $I$ to $S$. Now from the above definition we see that if $S$ corresponds to $G$, then $S^{(k)}_I$ corresponds to $G^{(k)}_v$, where $v = f^{-1}(I)$. This means that the operation sending $G$ to $G^{(k)}_v$ is closed on unit interval graphs, and hence they are expected to be $e$-positive. Motivated by Conjecture \ref{mconj}, we are led to study whether one can directly derive the $e$-positivity of $G^{(k)}_v$ from that of $G$ provided that $G$ is a unit interval graph. 

Similar issues occur for claw-free graphs. Gasharov (unpublished) and Stanley \cite{Sta98} conjectured that all claw-free graphs are $s$-positive. Note that if a graph $G$ is claw-free, then for any vertex $v$ of $G$ and any $k \ge 1$, $G^{(k)}_v$ is also claw-free. Inspired by Conjollary \ref{mconj-coro}, it is natural to consider whether one can directly derive the $s$-positivity of $G^{(k)}_v$ from that of $G$ provided that $G$ is claw-free.

Finally,  Conjecture \ref{mconj} might be valid for (claw, net)-free graphs. It is easy to see that if $G$ is (claw, net)-free, so is $G_v$ for any vertex $v \in V(G)$. As mentioned in the introduction, Foley, Ho\`{a}ng and Merkel \cite{FHM19} conjectured a graph is (claw, net)-free if and only if it is strongly $e$-positive. Provided the validity of this conjecture, it would be interesting if one can directly show that  the strong $e$-positivity of $G$ implies the $e$-positivity of $G_v$.

\medskip 

\noindent \textbf{Acknowledgments.} The third author is supported in part by the National Science Foundation of China (No. 11671037). The fourth author is supported in part by the Fundamental Research Funds for the Central Universities and the National Science Foundation of China (Nos. 11522110, 11971249).

\end{document}